\documentclass[12pt]{article}

\usepackage{amstext    }
\usepackage{amsthm    }
\usepackage{a4}
\usepackage[mathscr]{eucal}

\usepackage{amsmath}
\usepackage{amssymb}
\usepackage{amscd}

\numberwithin{equation}{section}

\newtheorem{theorem}{Theorem}[subsection]
\newtheorem{definition}[theorem]{Definition}
\newtheorem{proposition}[theorem]{Proposition}
\newtheorem{corollary}[theorem]{Corollary}
\newtheorem{lemma}[theorem]{Lemma}
\newtheorem{remark}[theorem]{Remark}

\newtheorem{example}[theorem]{Example}

\newcommand{\cali}[1]{\mathscr{#1}}

\newcommand{\const}{{\rm const}}
\newcommand{\dist}{{\rm dist}}

\newcommand{\loc}{{loc}}
\newcommand{\ddc}{dd^c}
\newcommand{\dc}{d^c}
\newcommand{\dbar}{\overline\partial}

\newcommand{\DSH}{{\rm DSH}}

\newcommand{\id}{{\rm id}}

\newcommand{\Cc}{\cali{C}}
\newcommand{\Dc}{\cali{D}}

\newcommand{\Fc}{\cali{F}}

\newcommand{\Kc}{\cali{K}}

\newcommand{\Mc}{\cali{M}}
\newcommand{\Nc}{\cali{N}}

\newcommand{\Rc}{\cali{R}}
\newcommand{\Uc}{\cali{U}}

\newcommand{\FS}{{\rm FS}}

\newcommand{\C}{\mathbb{C}}
\newcommand{\N}{\mathbb{N}}
\newcommand{\Z}{\mathbb{Z}}
\newcommand{\R}{\mathbb{R}}
\renewcommand{\P}{\mathbb{P}}

\renewcommand{\L}{{\cal L}}

\newcommand{\trans}{{\ \!\!\sp t\ \!\!}}

\title{Super-potentials for currents on compact K{\"a}hler manifolds
and dynamics of automorphisms}
\author{Tien-Cuong Dinh and Nessim Sibony}

\begin{document}
\maketitle

\begin{abstract}
We introduce a notion of super-potential (canonical function)
associated to positive closed $(p,p)$-currents on compact
K{\"a}hler manifolds and we develop a calculus on such currents. 
One of the key points in our study is the use of deformations in the space of
currents.
As an application, we obtain several results on the dynamics of holomorphic automorphisms:
regularity and uniqueness of the Green currents. We also get the regularity, the entropy, the ergodicity and the
hyperbolicity of the equilibrium measures.
\end{abstract}

\noindent
{\bf AMS classification :} 37F, 32H50, 32U40.

\noindent
{\bf Key-words :} super-potential, structural variety of currents,
moderate measure,
Green current, equilibrium measure.


\tableofcontents

\section{Introduction}

Let $(X,\omega)$ be a compact K{\"a}hler manifold of dimension $k$. Our
purpose in this paper is to develop a calculus on positive closed
currents of bidegree $(p,p)$ on $X$. We will also apply this calculus to
prove some surprising uniqueness results for dynamical currents
in their cohomology classes. 

When $S$ is a positive closed $(1,1)$-current on $X$, it is possible
to introduce its potential $u$ satisfying the following equation with
a normalization condition
$$\ddc u= S-\alpha \quad\mbox{and}\quad \int_X u\omega^k=0,$$
where $\alpha$ is a smooth representative of the cohomology class of
$S$. The function $u$ is quasi-p.s.h.; it is defined everywhere and satisfies $\ddc
u\geq -c\omega$ for some constant $c>0$. This is
the unique solution of the above equation and calculus problems on $S$ can be transfered
to computation on the potential $u$.

We have developed in \cite{DinhSibony10} a theory of super-potentials
associated to positive closed currents $S$ of bidegree $(p,p)$ in
$\P^k$ (our approach can be easily extended to homogeneous
manifolds). Let $\omega_\FS$ denote the Fubini-Study form on $\P^k$
normalized by $\int_{\P^k}\omega_\FS^k=1$. 
Assume for simplicity that $S$ is of mass 1, that is, $S$ is
cohomologous to $\omega_\FS^p$. One can always solve the
equation
\begin{equation}\label{eq_ddc}
\ddc U_S= S-\omega_\FS^p \quad\mbox{and}\quad \langle
U_S,\omega_\FS^{k-p+1}\rangle=0.
\end{equation}
But when $p>1$, the current $U_S$ is not unique and it is difficult to give
$U_S$ a value at every point, in order for example, to consider
expressions like the wedge-product $U_S\wedge [V]$ where $[V]$ is a current
associated to an analytic set $V$. In \cite{DinhSibony10}, we have
introduced for $S$ a super-potential
$\Uc_S$ which is a function defined on the space of positive closed
currents $R$ of bidegree $(k-p+1,k-p+1)$ and of mass 1. More
precisely, we have shown that it is possible to define\footnote{in
  \cite{DinhSibony10} we call $\Uc_S$ the super-potential of
  mean 0 of $S$.} 
$$\Uc_S(R)=\langle U_S,R\rangle:=\limsup \langle U_S,R'\rangle$$
with $R'$ smooth positive closed converging to $R$. 
The above formula is symmetric in $R$ and $S$, that is,
$\Uc_S(R)=\Uc_R(S)$. So, we also have $\Uc_S(R)=\langle S,U_R\rangle$
where $U_R$ is a normalized solution of the equation $\ddc U_R=R-\omega_\FS^{k-p+1}$. In particular, 
$\Uc_S(R)$ does not depend on the choice of $U_S$ and $U_R$.

The super-potentials appear as quasi-p.s.h. functions on an
infinite dimensional space and the value $-\infty$ is admissible. 
The calculus we have obtained is satisfactory and
permits to solve non trivial dynamical questions for holomorphic
endomorphisms of $\P^k$ and polynomial automorphisms of $\C^k$. It
also permits to give a useful intersection theory of positive closed
currents in $\P^k$. 

It will be important to extend such a calculus to arbitrary compact K{\"a}hler
manifolds. There are however some important difficulties. First,
according to Bost-Gillet-Soul{\'e} \cite{BostGilletSoule}, if $p>1$, it is not
always possible to solve the equation (\ref{eq_ddc}) with $U_S$ bounded
from above. In some sense, using the
potentials one may loose the positivity or the boundedness from below. Second, the approximation  of arbitrary positive
closed currents by smooth ones is only possible when a loss in
positivity is allowed, see Theorem \ref{th_regularization} below. The loss of
positivity is under control but it is still a source of several
technical difficulties. In general, the deformation of currents on
non-homogeneous manifolds is a delicate problem.

In the present paper, we introduce the super-potentials of $S$ as acting on
the real vector space of closed currents $R$ which are smooth and
cohomologous to 0. Then $\Uc_S$ is defined as
$$\Uc_S(R):=\langle S,U_R\rangle,$$
where $U_R$ is a smooth solution of $\ddc U_R=R$ satisfying 
some normalization conditions. This permits to develop the first steps
of a theory of super-potential on an arbitrary compact K{\"a}hler
manifold. In particular, we can define with some regularity
assumption, the wedge-product $S_1\wedge S_2$ where $S_j$ are positive
closed $(p_j,p_j)$-currents.

We then apply these notions to the dynamical study of automorphisms of
a compact K{\"a}hler manifold. Let $f$ be a holomorphic automorphism
of $X$. The dynamical degree of order $s$ of $f$ is defined as the
spectral radius of the pull-back operator $f^*$ on the cohomology
group $H^{s,s}(X,\R)$. It follows from a result by
Khovanskii-Teissier-Gromov \cite{Gromov1} that the function $s\mapsto
\log d_s$ is concave. In particular, we can assume that
$$1=d_0<d_1<\cdots< d_p=\cdots=d_{p'}>\cdots>d_{k-1}>d_k=1.$$
We have constructed in \cite{DinhSibony4} for $1\leq q\leq p$, Green
$(q,q)$-currents $T_q$. In our context, they are the positive closed currents of bidegree $(q,q)$ such that
$f^*(T_q)=d_qT_q$, see Section \ref{section_green} 
for the precise definition. Under the hypothesis that the dynamical degrees
are distinct (i.e. $p=p'$), we also constructed and studied an ergodic invariant
measure $\mu$ for $f$. 
The case of surfaces ($k=2$) was studied by Cantat in
\cite{Cantat}. Dynamically interesting automorphisms of surfaces are  also
constructed in Bedford-Kim \cite{BedfordKim} and McMullen \cite{McMullen}.

Here, we propose a new approach using super-potentials to deal with
convergence problems. We will show that the Green currents have
H{\"o}lder continuous super-potentials. The following uniqueness result
is quite surprising and can be applied to all $q$ smaller than or equal to
$p$. We refer to \cite{FornaessSibony3, FavreJonsson, Dujardin,
  DinhSibony9, DinhSibony10, DinhNguyenSibony1} and the references therein for analogous
results in other settings.

\medskip
\noindent
{\bf Theorem.} {\it Let $f$ be a holomorphic automorphism of a compact
  K{\"a}hler manifold $X$ and $d_s$ 
  the dynamical degrees of $f$. Suppose $V$ is a subspace of
  $H^{q,q}(X,\R)$  invariant under $f^*$. Assume that all the (real
  and complex) eigenvalues
  of the restriction of $f^*$ to $V$ are
  of  modulus strictly larger that $d_{q-1}$.
Then each class in $V$ contains at  most one positive closed
$(q,q)$-current.}

\medskip

As a consequence, we deduce that given a positive closed $(q,q)$-current $S$, the
convergence of the classes $(f^{n_i})^*[S]$, properly normalized,
implies the convergence of the currents $(f^{n_i})^*(S)$, properly
normalized. 
Here, $f^n:=f\circ\cdots\circ f$ ($n$ times) is the iterate of order
$n$ of $f$. The result applied to the current of integration on a
subvariety $Y$ of $X$ gives a description of the asymptotic behavior
of the inverse image of $Y$ by $f^n$ when $n\rightarrow\infty$.  
We also deduce that the Green currents are the unique positive closed currents in
their cohomology classes without restricting to invariant
currents. 

Assume that the dynamical degrees of $f$ are distinct,
i.e. $p=p'$ (for surfaces this just means $d_1>1$). Assume also  
that the action of $f^*$ on $H^{p,p}(X,\R)$ 
satisfies the following condition which is always true for surfaces. Let
$H$ be the invariant subspace of $H^{p,p}(X,\R)$ corresponding to
eigenvalues of maximal modulus. Suppose that $f^*$ restricted to
$H$ is diagonalizable over $\C$. This condition means that the Jordan form of $f^*$
restricted to $H\otimes_\R\C$ is a diagonal matrix. Let $T^+$ be a Green $(p,p)$-current of $f$ and
$T^-$ a Green $(k-p,k-p)$-current associated to $f^{-1}$. 
The hypothesis on  $f^*_{|H}$ is a
necessary and sufficient condition in order to have 
$T^+\wedge T^-\not=0$ for a suitable choice of $T^+$,
$T^-$, see Proposition \ref{prop_green_current_wedge}. These
measures $T^+\wedge T^-$ generate a real space $\Nc$ of finite dimension. We will
show that the convex cone $\Nc^+$ of positive measures in $\Nc$
is closed, with a simplicial basis and that the measures $\mu$ on the extremal
rays are ergodic. 
When the eigenvalues of $f^*_{|H}$ are all real positive, i.e. equal to $d_p$,
$\mu$ is mixing. 

We will show that any such measure $\mu$ is of maximal entropy $\log d_p$. 
Then, using a recent result of de Th{\'e}lin
\cite{deThelin} we deduce that $\mu$ is hyperbolic with precise
estimates on the positive and the negative Lyapounov exponents.  
The H{\"o}lder continuity of the super-potentials of the Green
currents implies that $\mu$ is moderate: if $u$ belongs to
a compact family of quasi-p.s.h. functions and $\ddc u\geq
-\omega$ then $\langle\mu, e^{\lambda|u|}\rangle\leq c$ for some
positive constants $\lambda$ and $c$. 
As far as we know, this property is the strongest regularity property
satisfied by the equilibrium measures in a quite general setting. It
implies that any quasi-p.s.h. function is in $L^p(\mu)$ for all $1\leq
p<\infty$. Moreover, $\mu$ has no mass on proper analytic subsets of
$X$. A result due to Katok \cite[p.694]{KatokHasselblatt}  
implies that the set of saddle periodic points is Zariski dense in $X$ since its  
 closure contains the support of $\mu$.

We have tried to make the paper readable for non experts. In Section
\ref{section_current} we give background on positive closed currents
and we introduce transforms on currents in a quite general context. We
show how to regularize positive closed currents and how to solve the
$\ddc$-equation in an arbitrary
compact K{\"a}hler manifold. In Section \ref{section_sp} we construct
an explicit
structural variety for a given current $R$ which is a difference of
positive closed currents. So, $R$ appears as the slice by $\{0\}\times X$
of a closed current $\Rc$ in $\P^1\times X$. This is the main technical tool,
which permits to use the powerful estimates on subharmonic functions
in order to prove the convergence theorems. We also define here the
intersection of currents. In the last section, we give the
applications to the dynamics of automorphisms. 

\medskip

\noindent
{\bf Main notations and conventions.} Throughout the paper, except for
some definitions 
in Section \ref{section_positive}, $(X,\omega)$ is a compact K{\"a}hler manifold of
dimension $k$. 
The notation $[V]$ or $[S]$ means the current of integration on an
analytic set $V$ or the class of a $\ddc$-closed $(p,p)$-current $S$ in
$H^{p,p}(X,\C)$ or $H^{p,p}(X,\R)$. 
Denote by 
$\pi:\widehat{X\times X}\rightarrow X\times X$ the blow-up along the diagonal
$\Delta$ and $\widehat\Delta:=\pi^{-1}(\Delta)$ the exceptional
hypersurface in $\widehat{X\times X}$.
The canonical projections of $X\times X$ on its factors
are denoted by $\pi_i$ and we define
$\Pi_i:=\pi_i\circ\pi$ for $i=1,2$. We also fix a K{\"a}hler form $\widehat\omega$
on $\widehat{X\times X}$. Let $\Cc_p$ denote the convex cone of positive closed
$(p,p)$-currents on $X$, $\Dc_p$ the real space generated by
$\Cc_p$ and $\Dc_p^0$ the subspace of currents in $\Dc_p$ which belong
to the class 0 in $H^{p,p}(X,\R)$. We consider on these spaces the
norms $\|\cdot\|_{\Cc^{-l}}$, $\|\cdot\|_*$ and the 
$\ast$-topology defined in Section \ref{section_positive}. On $\Cc_p$ or on
$\ast$-bounded (i.e. bounded with respect
to $\|\cdot\|_*$) subsets of $\Dc_p$, the $\ast$-topology coincides with the weak
topology on currents. The current $\Theta_0$, its deformations
$\Theta_\theta$ with $\theta\in\P^1$, the associated transforms $\L_0$,
$\L_\theta$ and the transform $\L_K$ are defined in Section \ref{section_regularization}.
The deformations $S_\theta:=\L_\theta(S)$ of a current $S$ and the associated
structural line $(S_\theta)_{\theta\in\P^1}$ are introduced in Sections \ref{section_regularization}
and \ref{section_structural}. The super-potential of a current $S$ in $\Dc_p$, normalized by
a fixed family $\alpha$ of closed $(p,p)$-forms, is denoted by $\Uc_S$. If
$S$ belongs to $\Dc_p^0$, then $\Uc_S$ does not depend on the choice
of $\alpha$, see Section \ref{section_sp_def}. Finally, most of the
constants depend only on $(X,\omega)$. The notations
$\gtrsim$, $\lesssim$ mean inequalities up to a multiplicative
constant and we will write $\sim$ when both inequalities are satisfied.


\section{Background on positive closed currents} \label{section_current}

In this section, we recall some basic facts on Hodge theory for
compact K{\"a}hler manifolds, some properties of positive closed
currents and plurisubharmonic functions. We refer
to \cite{Demailly3, deRham, Federer, Hormander, Lelong, Voisin} for
more detailed expositions on these subjects.

\subsection{Compact K{\"a}hler manifolds} \label{section_kahler}

\noindent
$\bullet$ {\bf Hodge cohomology groups.} Consider a compact K{\"a}hler manifold $X$ of
dimension $k$.
Let $H^r(X,\R)$ and $H^r(X,\C)$ denote the de Rham cohomology groups
of real and complex smooth $r$-forms. Let
$H^{p,q}(X,\C)$, $p+q=r$, be the subspace of $H^r(X,\C)$ generated by the classes of closed
$(p,q)$-forms. The Hodge theory asserts that
$$H^r(X,\C)=\bigoplus_{p+q=r} H^{p,q}(X,\C)\quad \mbox{and}\quad
H^{p,q}(X,\C)=\overline{H^{q,p}(X,\C)}.$$
For $p=q$, define
$$H^{p,p}(X,\R):=H^{p,p}(X,\C)\cap H^{2p}(X,\R),$$
then
$$H^{p,p}(X,\C)=H^{p,p}(X,\R)\otimes_\R\C.$$
The cup-product $\smile$ on $H^{p,p}(X,\R)\times
H^{k-p,k-p}(X,\R)$ is defined by
$$([\beta],[\beta'])\mapsto [\beta]\smile[\beta']:=\int_X\beta\wedge
\beta'$$
where $\beta$ and $\beta'$ are smooth closed forms. The last integral
depends only on the classes of $\beta$ and $\beta'$. The bilinear form $\smile$
is non-degenerate and induces a duality (Poincar\'e duality) between $H^{p,p}(X,\R)$ and
$H^{k-p,k-p}(X,\R)$. In the definition of $\smile$ one can take $\beta'$
smooth and $\beta$ a current in the sense of de Rham. So, 
$H^{p,p}(X,\R)$ can be defined as the quotient of the space of real closed
$(p,p)$-currents by the subspace of $d$-exact currents. Recall that a
$(p,p)$-current $\beta$ is real if $\overline \beta=\beta$. When
$\beta$ is a real $(p,p)$-current such that $\ddc\beta=0$, by the
$\ddc$-lemma \cite{Demailly3, Voisin}, the integral $\int_X\beta\wedge
\beta'$ is also independent of the choice of $\beta'$ smooth and
closed in a fixed cohomology class. So, using the duality, one can
associate to $\beta$ a class in $H^{p,p}(X,\R)$.

\bigskip

\noindent
$\bullet$ {\bf Blow-up along the diagonal.} 
The integration on the diagonal
$\Delta$ of $X\times X$ defines a real closed $(k,k)$-current $[\Delta]$ which is
positive, see
Section \ref{section_positive} for the notion of positivity.
By K{\"u}nneth formula, we have a canonical isomorphism
$$H^{k,k}(X\times X,\C)\simeq \sum_{0\leq r\leq k} H^{r,k-r}(X,\C)\otimes
H^{k-r,r}(X,\C).$$
Hence, $[\Delta]$ is cohomologous to a smooth real closed
$(k,k)$-form $\alpha_\Delta$ which is a finite combination of forms of
type $\beta(x)\wedge \beta'(y)$. Here, $\beta$ and $\beta'$ are closed
forms on $X$ of bidegree $(r,k-r)$ and $(k-r,r)$ respectively, and
$(x,y)$ denotes the coordinates of $X\times X$. 
In other words, if $\pi_i$ denote the projections of $X\times
X$ on its factors, then $\alpha_\Delta$ is a combination of
$\pi_1^*(\beta)\wedge\pi_2^*(\beta')$. So, $\alpha_\Delta$ satisfies 
$d_x\alpha_\Delta=d_y\alpha_\Delta=0$. 
Replacing $\alpha_\Delta(x,y)$ by ${1\over 2}\alpha_\Delta(x,y)+{1\over
  2}\alpha_\Delta(y,x)$ allows to assume that $\alpha_\Delta$ is
{\it symmetric}, i.e. invariant by the involution $(x,y)\mapsto (y,x)$.

Let $\pi:\widehat{X\times X}\rightarrow X\times X$ be the blow-up of
$X\times X$ along $\Delta$ and $\widehat\Delta:=\pi^{-1}(\Delta)$ the
exceptional hypersurface. By a theorem of Blanchard
\cite{Blanchard}, $\widehat{X\times X}$ is a compact K{\"a}hler
manifold. We fix a K{\"a}hler form $\widehat\omega$ on $\widehat{X\times X}$. 
According to Gillet-Soul{\'e} 
\cite[1.3.6]{GilletSoule}, there is a real smooth closed $(k-1,k-1)$-form
$\eta$ on $\widehat{X\times X}$ such that $\pi^*(\alpha_\Delta)$ is cohomologous to
$[\widehat\Delta]\wedge\eta$, where $[\widehat\Delta]$ is the positive
closed $(1,1)$-current
of integration on $\widehat\Delta$. Hence, $\pi_*([\widehat\Delta]\wedge\eta)$
is cohomologous to $\alpha_\Delta$ and to $[\Delta]$. On the other hand, 
$\pi_*([\widehat\Delta]\wedge\eta)$ is supported on $\Delta$ and is equal to a product of
$[\Delta]$ by a function. We deduce that $\pi_*([\widehat\Delta]\wedge\eta)=[\Delta]$.
The map $(x,y)\mapsto (y,x)$ induces an involution on
$\widehat{X\times X}$. We can also choose $\eta$ symmetric with
respect to this involution.

Let $\gamma$ be a real closed $(1,1)$-form on $\widehat{X\times X}$,
cohomologous to $[\widehat\Delta]$. We can choose $\gamma$
symmetric. We will see later that there is a 
quasi-p.s.h. function $\varphi$ on $\widehat{X\times X}$ such that
$\ddc\varphi=[\widehat\Delta]-\gamma$. This function is necessarily
symmetric. Subtracting from $\varphi$ a
constant allows to assume that $\varphi<-2$.

\medskip

\noindent
$\bullet$ {\bf Local coordinates near $\Delta$ and
$\widehat\Delta$.} Consider a local coordinate system $x=(x_1,\ldots,x_k)$
on a chart of $X$. For simplicity assume that the ball $W$ of center 0 and
of radius 1 is strictly contained in this chart.
For the neighbourhood $W\times W$ of $(0,0)$ in $X\times X$, we will
use the coordinates $(x,y)=(x_1,\ldots,x_k,y_1,\ldots,y_k)$. 
The diagonal $\Delta$ contains the
point $(0,0)$ and is given by the equation $x=y$. Define
$x':=x-y$. Then $(x',y)$ is also a coordinate system of $W\times
W$ and $\Delta$ is given by $x'=0$. Consider the
submanifold $M$ of $\C^k\times\C^k\times\P^{k-1}$ defined by
$$M:=\big\{(x',y,[v])\in \C^k\times\C^k\times\P^{k-1},\quad x'\in [v]\big\},$$
where $[v]=[v_1:\cdots:v_k]$ denotes the homogeneous coordinates of
$\P^{k-1}$. Recall that $x'$ belongs to $[v]$ if and only if $x'$ and
$v$ are proportional. The submanifold $M$ is the blow-up of
$\C^k\times\C^k$ along $x'=0$. So, we identify $\pi^{-1}(W\times W)$
with an open set in $M$ defined by $\|x'+y\|=\|x\|<1$ and $\|y\|<1$.

Consider a point $(a,b,[u])$ in $\widehat\Delta$. We have necessarily
$a=0$. For simplicity, assume that the first coordinate of $u$ is the
largest one.
Therefore, we can write $[u]=[1:u_2:\cdots:u_k]$ with $|u_i|\leq 1$. In a
neighbourhood of $(0,b,[u])$, the first coordinate of $v$ does not
vanish and we can write $[v]=[1:v_2:\cdots:v_k]$ with $|v_i|<2$. Write
$v':=(v_2,\ldots,v_k)$. Then,
$(x_1',y,v')$ is a local coordinate system for a
neighbourhood of $(0,b,[u])$.  Here, $\widehat\Delta$ is given by the
equation $x'_1=0$. We also have
$$\pi(x'_1,y,v')=(x',y)=(x'_1,x'_1v',y)\quad
\mbox{and}\quad \Pi_2(x'_1,y,v')=y.$$
We see that $\Pi_2$ and its restriction to $\widehat\Delta$ are submersions. In the same way, we prove that
$\Pi_1$ and its restriction to $\widehat\Delta$ are also submersions.

\subsection{Positive currents and plurisubharmonic functions} \label{section_positive}

\noindent
$\bullet$ {\bf Positive closed currents.}
A smooth $(p,p)$-form $\phi$ on a general complex manifold
of dimension $k$ is
{\it positive} if it can be written in local charts as a finite combination
with positive coefficients of forms of type 
$$(i\alpha_1\wedge \overline\alpha_1)\wedge \ldots\wedge
(i\alpha_p\wedge \overline\alpha_p)$$
where $\alpha_i$ are $(1,0)$-forms. The positivity is a pointwise
property and does not depend on local coordinates. A
$(p,p)$-current $S$ is {\it weakly positive} if $S\wedge\phi$ is a
positive measure for every smooth positive $(k-p,k-p)$-form
$\phi$. The current $S$ is {\it positive} if $S\wedge\phi$ is a
positive measure for every smooth weakly positive $(k-p,k-p)$-form
$\phi$. The notions of positivity and weak positivity coincide for
$p=0,1,k-1$ and $k$. We say that $S$ is {\it negative} if $-S$ is
positive and we write $S\geq S'$, $S'\leq S$ when $S-S'$ is positive.
Note that positive and negative currents are real. If $S$, $S'$ are
positive and $S'$ is smooth then $S\wedge S'$ is positive.
Let $V$ be an analytic subset of pure codimension $p$. Then, the
integration on the regular part of $V$ defines a positive closed
$(p,p)$-current that we denote by $[V]$. A $(p,p)$-current $S$ is said to be
{\it strictly positive} if in local coordinates $x$, we have $S\geq
\epsilon (\ddc\|x\|^2)^p$ for some $\epsilon>0$. 

Let $(X,\omega)$ be a compact K{\"a}hler manifold of
dimension $k$.
If $S$ is a positive or negative $(p,p)$-current on $X$, define {\it the
  mass}\footnote{in this case, this quantity is equivalent to the mass norm for
  currents of order 0, see \cite{Federer}.} of $S$ by
$$\|S\|:=|\langle S,\omega^{k-p}\rangle|.$$
Let $\Cc_p$ denote the cone of positive closed $(p,p)$-currents on
$X$, $\Dc_p$ the real space generated by $\Cc_p$ and $\Dc^0_p$
the space of currents $S\in \Dc_p$ such that $[S]=0$ in $H^{p,p}(X,\R)$.
The duality between the cohomology groups implies that 
if $S$ is a current in $\Cc_p$, its mass depends only
on the class $[S]$ in $H^{p,p}(X,\R)$. Define {\it the norm
  $\|\cdot\|_*$} on $\Dc_p$ by
$$\|S\|_*:=\min\|S^+\|+\|S^-\|,$$
where the minimum is taken over $S^+$, $S^-$ in $\Cc_p$ such that
$S=S^+-S^-$. 
A subset in $\Dc_p$ is {\it $\ast$-bounded} if it is bounded with
respect to the $\|\cdot\|_*$-norm. 
We will consider on $\Dc_p$ and $\Dc^0_p$ the following
 {\it $\ast$-topology}. We say that {\it $S_n$ converge to $S$
in $\Dc_p$} if $S_n\rightarrow S$ weakly and if $\|S_n\|_*$ is
bounded by a constant independent of $n$. Note that the $\ast$-topology
restricted to $\Cc_p$ or to a $\ast$-bounded subset of $\Dc_p$ coincides with the weak topology.
We will see in
Theorem \ref{th_regularization}  
below that smooth forms are dense in $\Dc_p$ and $\Dc^0_p$ for the $\ast$-topology.

Consider some natural {\it norms} and {\it distances} on $\Dc_p$. 
For $l\geq 0$, let $[l]$ denote the integer part of $l$. Let
$\Cc_{p,q}^l$ be the space of $(p,q)$-forms whose coefficients admit all
derivatives of order $\leq [l]$ and these derivatives are $(l-[l])$-H{\"o}lder continuous.  
We use here the sum of $\Cc^l$-norms of the coefficients for a
fixed atlas.
If $S$ and $S'$ are
currents in $\Dc_p$, define\footnote{the definition is meaningful for
  any current $S$ of order 0 and $\|\cdot\|_{\Cc^{-0}}$ is equivalent
  to the mass norm in the usual
  sense, see \cite{Federer}.} 
$$\|S\|_{\Cc^{-l}}:= \sup_{\|\Phi\|_{\Cc^l}\leq 1} |\langle
S,\Phi\rangle| \quad \mbox{and}\quad 
\dist_l(S,S'):=\|S-S'\|_{\Cc^{-l}}$$
where $\Phi$ is a test smooth $(k-p,k-p)$-form on $X$. 
Observe that $\|\cdot\|_{\Cc^{-l}}\lesssim \|\cdot\|_*$ for every
$l\geq 0$. 
The following result
is proved as in \cite{DinhSibony10} using the theory of interpolation
between Banach spaces.

\begin{proposition} \label{prop_compare_dist}
Let $l$ and $l'$ be real strictly positive numbers with $l<l'$. 
Then on any $\ast$-bounded subset of $\Dc_p$, the topology  induced by 
$\dist_l$ or by $\dist_{l'}$ coincides with the weak topology.
Moreover, on any $\ast$-bounded subset of $\Dc_p$, there is a constant
$c_{l,l'}>0$ such that 
$$\dist_{l'}\leq\dist_l\leq
c_{l,l'}[\dist_{l'}]^{l/l'}.$$
In particular, a function on a $\ast$-bounded subset of $\Dc_p$ is H{\"o}lder continuous with
respect to $\dist_l$ if and only if it is H{\"o}lder continuous with respect to $\dist_{l'}$.
\end{proposition}

\noindent
$\bullet$ {\bf Plurisubharmonic functions.}
Consider a general (connected) complex manifold $X$. An upper semi-continuous function $u:X\rightarrow
\R\cup\{-\infty\}$, not identically $-\infty$, is {\it
plurisubharmonic} (p.s.h. for short) if its restriction to each
holomorphic disc in $X$ is subharmonic or identically equal to
$-\infty$. If $u$ is a p.s.h. function then $u$ belongs to $L^p_\loc$
for $1\leq p<\infty$, and $\ddc u$ is a
positive closed $(1,1)$-current on $X$. 
Conversely, if $S$ is a positive closed $(1,1)$-current, it can be
locally written as $S=\ddc u$ with $u$ p.s.h.
A subset of $X$ is {\it locally pluripolar} if it is locally
contained in the pole set $\{u=-\infty\}$ of a p.s.h. function.
P.s.h. functions satisfy a
maximum principle. In particular, on a compact manifold,
p.s.h. functions are constant. A function $u$ on $X$ is {\it
  quasi-p.s.h.} if it is locally a difference of a p.s.h. function and
a smooth function.

Assume now that $X$ is a compact K{\"a}hler manifold of dimension $k$ and $\omega$ is a
K{\"a}hler form on $X$. If $u$ is a quasi-p.s.h. function on $X$ then  $\ddc u+c\omega$ is a
positive closed $(1,1)$-current for $c>0$ large enough. Conversely, if
$S$ is a positive closed $(1,1)$-current and $\alpha$ is a real closed
smooth $(1,1)$-form cohomologous to $S$, then there is a quasi-p.s.h. function $u$ such
that $\ddc u=S-\alpha$. 
The function $u$ is unique up to an additive constant.
A subset of $X$ is {\it
  pluripolar} if it is contained in the pole set $\{u=-\infty\}$ of a
quasi-p.s.h. function $u$. 

A function is called {\it d.s.h.} if it is
equal outside a pluripolar set to a difference of two
quasi-p.s.h. functions. We identify two d.s.h. functions if they are
equal out of a pluripolar set. If $u$ is d.s.h. then $\ddc u$ is a
difference of two positive closed $(1,1)$-currents which are cohomologous. Conversely, if
$S^+$ and $S^-$ are positive closed $(1,1)$-currents in the same
cohomology class then $S^+-S^-=\ddc u$ for some d.s.h. function
$u$. The function $u$ is unique up to an additive constant.
There are several equivalent norms on the space of
d.s.h. functions. We consider the following one, see \cite{DinhSibony6}
$$\|u\|_\DSH:=\|u\|_{L^1}+\|\ddc u\|_*.$$
We have the following proposition \cite{DinhSibony6}.

\begin{proposition} Let $u$ be a d.s.h. function on $X$. Then there
  exist two quasi-p.s.h. functions $u^+$, $u^-$ such that
$$u=u^+-u^-,\quad \|u^\pm\|_{L^1}\leq c\|u\|_\DSH,\quad\mbox{and}\quad \ddc u^\pm\geq -c\|u\|_\DSH\omega,$$
where $c>0$ is a constant independent of $u$.
\end{proposition}

We deduce from this proposition and the fundamental exponential estimate for
p.s.h. functions \cite{Hormander} the following result, see also \cite{DinhNguyenSibony}.

\begin{proposition} \label{prop_exp}
 There are constants $\lambda>0$ and $c>0$ such
  that if $u$ is a d.s.h. function with $\|u\|_\DSH\leq 1$ then
$$\int_X e^{\lambda|u|}\omega^k\leq c.$$ 
\end{proposition}

We will need the following version of the exponential estimate for
d.s.h. functions on $\P^1$ and for $\omega_\FS$ the
Fubini-Study form on $\P^1$.

\begin{lemma} \label{lemma_dsh_exp}
Let $u$ be a d.s.h. function on $\P^1=\C\cup\{\infty\}$.
  Assume that $u$ vanishes outside
  the unit disc of $\C$ and that $\ddc u$ is a measure of mass
  at most equal to $1$. Then there are constants $\lambda>0$ and $c>0$
  independent of $u$ such that
$$\int_{\P^1} e^{\lambda|u|}\omega_\FS\leq c,$$
In particular, if $B$ is a disc of radius $r$, $0<r<1/2$, then $\inf_B
|u|\leq -c'\log|r|$ for some constant $c'>0$ independent of $u$, $B$
and $r$. 
\end{lemma}
\proof
Write $\ddc u=\nu^+-\nu^-$ where $\nu^\pm$ are probability measures
with support in the unit disc. Define for $z\in\C$
$$u^\pm(z):=\int_\C\log|z-\xi|d\nu^\pm(\xi).$$ 
Observe that $\|u^\pm\|_{L^1(\P^1)}$ are bounded by a constant independent of
$\nu^\pm$. We also have 
$$\lim_{z\rightarrow\infty}
  u^\pm(z)-\log|z|=0\quad \mbox{and}\quad \ddc
u^\pm=\nu^\pm-\delta_\infty$$ 
where $\delta_\infty$ is the Dirac mass
at $\infty$. It follows that 
$$\lim_{z\rightarrow\infty} u^+(z)-u^-(z)=0\quad\mbox{and}\quad
\ddc(u^+-u^-)=\nu^+-\nu^-=\ddc u.$$
So, $u^+-u^-$ and $u$ differe by a constant. The fact that $u$ is
supported in the unit disc implies that $u=u^+-u^-$. We deduce that
$\|u\|_{L^1}$ is bounded by a constant independent of $u$, and then
$\|u\|_\DSH$ is bounded by a constant independent of $u$. Proposition
\ref{prop_exp} implies the result.
\endproof

\noindent
$\bullet$ {\bf Slicing of positive closed current.}
Let $V$ be a complex manifold of dimension $l$. We are interested in families of currents
parametrized by $V$ which are slices of some closed current $\Rc$ in
$V\times X$. Let $\pi_V$
and $\pi_X$ denote the canonical projections from $V\times X$ on its
factors. We have the
following proposition where currents on $\{\theta\}\times X$
are identified with currents on $X$, see also \cite{DinhSibony7}.

\begin{proposition} \label{prop_slice}
Let $\Rc$ be a positive closed $(s,s)$-current in $V\times
X$ with $s\leq k$. Then there is a locally pluripolar subset $E$ of $V$ such that the slice
$\langle \Rc,\pi_V,\theta\rangle$ exists for $\theta\in V\setminus
E$. Moreover, $\langle \Rc,\pi_V,\theta\rangle$ is a positive closed
$(s,s)$-current on $\{\theta\}\times X$ and its class in
$H^{s,s}(X,\R)$ does not depend on $\theta$. 
\end{proposition}

Recall that slicing is the generalization of restriction of forms to level sets of
holomorphic maps. It can be viewed as a version of Fubini's theorem or
Sard's theorem
for currents. The operation is well-defined for currents $\Rc$ of order 0 and of
bidegree $\leq (k,k)$ such that
$\partial \Rc$ and $\dbar\Rc$ are of order 0. When $\Rc$ is a smooth form, 
$\langle \Rc, \pi_V,\theta\rangle$ is simply the restriction of $\Rc$
to $\{\theta\}\times X$. 
When $\Rc$ is the current of integration on an analytic subset $Y$ of
$V\times X$, 
$\langle \Rc,\pi_V,\theta\rangle$ is the current of integration on the analytic set 
$Y\cap \{\theta\}\times X$ for $\theta$ generic. 

In general, if $\phi$ is a
smooth form on $V\times X$ then
$\langle \Rc\wedge\phi,\pi_V,\theta\rangle=\langle \Rc,\pi_V,\theta\rangle\wedge\phi$.     
Slicing commutes with the operations $\partial$
and $\overline\partial$. So, in our situation, since $\Rc$ is closed, $\langle \Rc,\pi_V,\theta
\rangle$ is also closed. The following description shows that $\langle \Rc,\pi_V,\theta
\rangle$ is positive.

Let $z$ denote the coordinates in a chart of $V$ and $\lambda_V$ the standard volume form.
Let $\psi(z)$ be a positive smooth function with compact support
such that $\int\psi\lambda_V=1$. Define
$\psi_\epsilon(z):=\epsilon^{-2l}\psi(\epsilon^{-1} z)$ and 
$\psi_{\theta,\epsilon}(z):=\psi_\epsilon(z-\theta)$. The measures $\psi_{\theta,\epsilon}\lambda_V$
approximate the Dirac mass at $\theta$. For every smooth test form $\Psi$
of bidegree $(k-s,k-s)$ on $V\times X$
one has
\begin{equation} \label{eq_slice}
\langle \Rc,\pi_V,\theta\rangle (\Psi)=\lim_{\epsilon\rightarrow 0}
\langle \Rc\wedge \pi_V^*(\psi_{\theta,\epsilon}\lambda_V),\Psi\rangle
\end{equation}
when $\langle \Rc,\pi_V,\theta\rangle$ exists.
This property holds for all choice of $\psi$.
Conversely, when the previous limit exists and is independent of $\psi$, 
it defines the current $\langle \Rc,\pi_V,\theta\rangle$ and one says 
that $\langle \Rc,\pi_V,\theta\rangle$ {\it is well-defined}.
The following formula holds for smooth forms $\Omega$ 
of maximal degree with compact support in $V$:
\begin{eqnarray} \label{eq_slice_bis}
\int_{\theta\in V}\langle \Rc,\pi_V,\theta\rangle (\Psi)\Omega(\theta) & = &  \langle \Rc\wedge \pi_V^*(\Omega),
\Psi\rangle.
\end{eqnarray}

\noindent
{\bf Proof of Proposition \ref{prop_slice}.} Since the problem is
local on $V$,
we can assume that $V$ is a ball in $\C^l$ and $z$ are the standard
coordinates. It is enough to consider real test forms $\Psi$ with
compact suppport. Define 
$\phi:=(\pi_V)_*(\Rc\wedge \Psi)$. This is a current of bidegree
$(0,0)$ on $V$.  Observe
that $\ddc \Psi$ can be written as a difference of positive closed
forms on $V\times X$, not necessarily with compact support. It follows
that $\ddc \phi=(\pi_V)_*(\Rc\wedge\ddc\Psi)$ is a difference of positive closed
currents. Therefore, $\phi$ can be considered as a d.s.h. function. We have
$$\langle \Rc\wedge
\pi_V^*(\psi_{\theta,\epsilon}\lambda_V),\Psi\rangle = \int_V
\phi\psi_{\theta,\epsilon}\lambda_V.$$
Classical properties of p.s.h. functions imply that for $\theta$
outside a pluripolar set (because $\phi$ is only d.s.h.) the last integral
converges to $\phi(\theta)$ 
when $\epsilon\rightarrow 0$. So, for such a $\theta$, 
the limit in (\ref{eq_slice}) exists
and does not depend on $\psi$. 

Choose a pluripolar set $E\subset V$
such that the previous convergence holds for $\theta$ outside $E$ and
for a countable family $\Fc$ of test forms $\Psi$. We choose a family $\Fc$ which is dense for the
$\Cc^0$-topology. The density implies that
we have the convergence for a test form $\Psi$ strictly
positive near $\{\theta\}\times X$. This,
the density of $\Fc$ together with the positivity of $\Rc$ imply the convergence for every $\Psi$. 
Hence, $\langle\Rc,\pi_V,\theta\rangle$ is well-defined
 for $\theta\not\in E$ and is a positive closed current on
 $\{\theta\}\times X$. 

We have $\langle \Rc,\pi_V,\theta\rangle (\Psi)=\phi(\theta)$. 
Consider a closed $(k-s,k-s)$-form $\Phi$ on $X$ and
$\Psi:=\pi_X^*(\Phi)$. If $\phi$ is defined as above, we have
$d\phi=0$. Therefore, $\phi$ is a constant function and 
$\langle\Rc,\pi_V,\theta\rangle(\Psi)$ does not
depend on $\theta$. If $\langle \Rc,\pi_V,\theta\rangle$ is identified
with a current on $X$, then $\langle \Rc,\pi_V,\theta\rangle(\Phi)$ is
independent of $\theta$. This, together with Poincar{\'e}'s duality,
implies that the class of
$\langle \Rc,\pi_V,\theta\rangle$ in $H^{s,s}(X,\R)$ does not depend
on $\theta$.
\hfill $\square$

\begin{remark} \label{rk_slice_continuity}
\rm
Assume that $\langle \Rc,\pi_V,\theta\rangle$ is defined for $\theta$ outside a set
of zero measure and that $\theta\mapsto \langle
\Rc,\pi_V,\theta\rangle$ can be extended to a continuous map with
values in the space of currents of order 0. 
Then, by definition of slicing, (\ref{eq_slice}) and (\ref{eq_slice_bis}), $\langle
\Rc,\pi_V,\theta\rangle$ is defined for every $\theta$ and coincides
with the continuous extension of  $\theta\mapsto \langle
\Rc,\pi_V,\theta\rangle$. If $\Rc_n$ are positive closed currents
converging to $\Rc$, we can prove that there is a subsequence
$\Rc_{n_i}$ with $\langle \Rc_{n_i},\pi_V,\theta\rangle\rightarrow
\langle \Rc,\pi_V,\theta\rangle$ for almost every $\theta$. Indeed,
for a bounded sequence of d.s.h. functions on $V$ we can extract a
subsequence which converges almost everywhere. 
\end{remark}

\subsection{Transforms on currents} \label{section_transform}

$\bullet$ {\bf General transforms on currents.} 
We recall here a general idea how to construct linear operators
on currents which are useful in geometrical questions. 
Let $X_1$, $X_2$ and $Z$ be Riemannian manifolds and $\tau_1$, $\tau_2$
smooth maps from $Z$ to $X_1$ and $X_2$. Let $\Theta$ be a fixed
current on $Z$. Define for a current $S$ on $X_1$ another current
$\L_\Theta(S)$ on $X_2$ by 
$$\L_\Theta(S):=(\tau_2)_*\big(\tau_1^*(S)\wedge\Theta\big)$$
when the last expression is meaningful. The current $\tau_1^*(S)$ is
well-defined if $S$ is a bounded form or if $\tau_1$ is a
submersion. The operator $(\tau_2)_*$ is well-defined if $\tau_2$ is
proper, in particular, when $Z$ is compact. Assume that $Z$ is
compact. Then, $\L_\Theta$ is well-defined on smooth currents $S$. Let
$\Theta'$ denote the push-forward of $\Theta$ to $X_1\times X_2$ by
the map $(\tau_1,\tau_2)$. Then, $\Theta'$ defines a transform
$\L_{\Theta'}$ where $Z$ is replaced by $X_1\times X_2$ . The
transform $\L_{\Theta'}$ 
is equal to $\L_\Theta$ on smooth forms $S$ and this useful property may be
extended to larger spaces of currents.

In this paper, we consider the following situation used in
Gillet-Soul{\'e} \cite{GilletSoule} and \cite{BostGilletSoule,
  DinhSibony3, DinhSibony4, Vigny}, see also \cite{Berndtsson, Henkin, Meo2, DinhSibony6}.
We use the notations introduced in Section \ref{section_kahler}.
Consider
a current $\Theta$ of bidegree $(r,s)$ on $\widehat {X\times X}$. If $S$ is a
current on $X$, define the transform $\L_\Theta(S)$ of $S$
by 
$$\L_{\Theta}(S):=(\Pi_2)_*\big(\Pi_1^*(S)\wedge \Theta\big).$$
This definition makes sense if the last wedge-product is well-defined, in particular when
$\Theta$ or $S$ is smooth.  

If $S$ is of bidegree $(p,q)$ then
$\L_{\Theta}(S)$ is of bidegree $(p+r-k,q+s-k)$. So, we say
  that the transform $\L_{\Theta}$ is {\it of bidegree }
  $(r-k,s-k)$. The bidegree may be negative. In what follows, we will
  be interested in 
  the cases where $r=s=k$ or $r=s=k-1$, and $S$ is a
  current in $\Dc_p$. The current $\Theta$ will be real and smooth or
  smooth outside  $\widehat\Delta$.   If $\Theta$ is positive or
  negative, we say that the transform $\L_\Theta$ is {\it positive} or
  {\it negative} respectively. 

\begin{example} \label{example_id} \rm
Consider $\Theta_0:=[\widehat\Delta]\wedge\eta$ where $\eta$ is the
smooth real closed form of bidegree $(k-1,k-1)$ chosen in Section
\ref{section_kahler}, and define $\L_0:=\L_{\Theta_0}$. 
Since $\Pi_1$ and its restriction to $\widehat\Delta$ are submersions,
$\L_0$ can be extended continuously to any current $S$.
We have  $\pi_*(\Theta_0)=[\Delta]$.
So, if $S$
is a smooth form, then
$$\L_0(S)=(\pi_2)_*\big(\pi_1^*(S)\wedge [\Delta]\big)=S.$$
By continuity, $\L_0$ is equal to the identity on all currents $S$. 
If $S$ is in $\Dc_p$, using the theory of intersection with positive
closed $(1,1)$-currents \cite{Demailly2, Demailly1, FornaessSibony2}
and that $[\widehat\Delta]=\ddc\varphi+\gamma$, we obtain for $S$ in
$\Dc_p$ that
if $S':=\Pi_1^*(S)$
$$S=\L_0(S)=(\Pi_2)_*\big(\ddc(\varphi S'\wedge \eta)+\gamma\wedge
S'\wedge \eta\big),$$
see also \cite{DinhSibony3}. 
We
will construct some deformations $\L_\theta$ of $\L_0$,
i.e. transforms associated to some deformations $\Theta_\theta$ of $\Theta_0$.
\end{example}

\noindent
$\bullet$ {\bf Regular and semi-regular transforms.}
Consider now a situation used in \cite{DinhSibony3, DinhSibony4}. 
Let $\Theta$ be a form which is smooth outside $\widehat\Delta$ and such that
$$|\Theta|\lesssim -\log\dist(\cdot,\widehat\Delta)\quad
\mbox{and} \quad 
|\nabla\Theta|\lesssim \dist(\cdot,\widehat\Delta)^{-1}$$ 
near $\widehat\Delta$. 
Here, the estimate on $\nabla\Theta$ means an estimate on the
gradients of the coefficients of $\Theta$ for a fixed atlas.
Transforms associated to such forms $\Theta$ are called {\it
  semi-regular} (when $\Theta$ is smooth, we say that
  $\L_\Theta$ is {\it regular}). 
The form $\Theta':=\pi_*(\Theta)$ is 
smooth outside $\Delta$. 
Using the local coordinates described in Section \ref{section_kahler}, one
proves that 
$$|\Theta'|\lesssim
-\log\dist(\cdot,\Delta)\dist(\cdot,\Delta)^{2-2k}$$
and
$$|\nabla\Theta'|\lesssim
\dist(\cdot,\Delta)^{1-2k}$$
near $\Delta$, see \cite{DinhSibony3}.
In particular, the coefficients of $\Theta'$ restricted to
$X\times\{y\}$ are in $L^{1+1/k}$ for every $y\in X$. 

Recall that $\L_\Theta$ is defined for $S$ smooth. Since $\Pi_1^*(S)\wedge \Theta$ has no mass on
$\widehat\Delta$, we have
$$\L_\Theta(S)=(\pi_2)_*\big(\pi_1^*(S)\wedge \Theta').$$
The wedge-product $\pi_1^*(S)\wedge\Theta'$ has no mass on
$\Delta$. So, one has to integrate only outside $\Delta$. Then,
using the estimates on
$|\Theta'|$, $|\nabla\Theta'|$ and the H{\"o}lder inequality, we obtain the following
result, see \cite{DinhSibony3, DinhSibony4}.

\begin{proposition} \label{prop_regular}
Any semi-regular transform can be extended to a linear
continuous operator from the space of
currents of order $0$ to the space of $L^{1+1/k}$ forms.
It defines a linear continuous operator from the space of $L^q$
forms, $q\geq 1$, to the space of $L^{q'}$ forms where $q'$ is
given by $q'^{-1}+1=q^{-1}+[1+1/k]^{-1}$ if $q< k+1$ and
$q'=\infty$ if $q\geq k+1$. It also defines a linear continuous
operator from the space of $L^\infty$ forms to the space of $\Cc^1$ forms.
\end{proposition}

The following result is a direct consequence of Proposition \ref{prop_regular}.

\begin{corollary} \label{cor_regular}
Let $\L_1$, $\ldots$, $\L_{k+2}$ be semi-regular transforms of bidegree
$(0,0)$. If $S$ is a current of order $0$, then
$S':=\L_{k+2}\circ\cdots\circ \L_1(S)$ is a form of class
$\Cc^1$. Moreover, we have
$\|S'\|_{\Cc^1}\leq c\|S\|$, where $c>0$ is a constant independent of
$S$. 
\end{corollary}

We will need the following lemma.

\begin{lemma} \label{lemma_tr_norm}
Assume that $\Theta$ is a smooth positive closed $(k,k)$-form.
Then $\L_\Theta$ defines a linear map from $\Dc_p$ to itself which preserves
$\Cc_p$, $\Dc^0_p$ and is continuous with respect to the $\ast$-topology. 
Moreover, if $S$ is in $\Dc_p$, then $\|\L_\Theta(S)\|_*\leq c \|\Theta\|\|S\|_*$
for some constant $c>0$ independent of $\Theta$ and $S$.
\end{lemma}
\proof
Since $\Pi_1$ is a submersion, by definition, $\L_\Theta$ is a linear
continuous map on currents. It is clear that $\L_\Theta$ preserves
$\Cc_p$, $\Dc_p$ and $\Dc^0_p$. We only have to prove the estimate on $\|\L_\Theta(S)\|_*$.
We can assume that 
 $S$ is positive. Since $\Pi_1$ is a
submersion, we have $\|\Pi_1^*(S)\|\lesssim \|S\|$.
Recall that the mass of a positive closed current
can be computed cohomologically. Therefore, 
$$\|\Pi_1^*(S)\wedge\Theta\|\lesssim \|\Theta\|\|\Pi_1^*(S)\|\lesssim
\|\Theta\|\|S\|.$$
The continuity of $(\Pi_2)_*$ implies the result. 
\endproof

\noindent
$\bullet$ {\bf Symmetric transforms.} The map $(x,y)\mapsto
(y,x)$ on $X\times X$ induces an involution on $\widehat{X\times X}$. In order to
simplify notations, we will only consider the transforms $\L_\Theta$ associated to
forms $\Theta$ which are invariant by this involution. We say that
$\L_\Theta$ is {\it symmetric}. Let $\Psi$ be a smooth test form on $X$ of the
apropriate bidegree. If $S$ is smooth then we deduce from the symmetry of
$\L_\Theta$ that
$$\langle \L_\Theta(S),\Psi\rangle =\int_{\widehat{X\times X}}
\Pi_1^*(S)\wedge \Theta\wedge\Pi_2^*(\Psi) = \langle
S, (\Pi_1)_*(\Pi_2^*(\Psi)\wedge\Theta)\rangle = \langle S, \L_\Theta(\Psi)\rangle.$$
When $S$
is closed and $\Phi$ is a smooth test form of the apropriate bidegree, we have
$$\langle \L_\Theta(S),\ddc \Phi\rangle = \int_{\widehat{X\times X}}
\Pi_1^*(S)\wedge \ddc\Theta\wedge\Pi_2^*(\Phi)=\langle \L_{\ddc\Theta}(S),\Phi\rangle.$$
Observe that the smoothness of $S$ is superflous when $\Theta$ is
smooth. The previous identities may be extended to some cases where
$S$ and $\Theta$ are not smooth using a regularization on $S$. 

\subsection{Regularization and Green potential} \label{section_regularization}

\noindent
$\bullet$ {\bf Deformation of the identity transform.}
We introduce in this section a family of regular transforms
$\L_\theta$, $\theta\in\P^1\setminus\{0\}$, of bidegree $(0,0)$
which is a continuous deformation of the identity transform $\L_0$
considered in Example \ref{example_id}.

We use the notations of
Section \ref{section_kahler}.
Consider the following
regularization of the function $\varphi$. Recall that $\varphi\leq -2$. Let $\chi$ be a smooth convex
increasing function on $\R\cup\{-\infty\}$ such that $\chi(t)=t$ for $t\geq 0$,
$\chi(t)=-1$ for $t\leq -2$ and $0\leq\chi'\leq 1$. Define
$$\chi_\theta(t):=\chi(t-\log|\theta|)+\log|\theta|\quad\mbox{and}
\quad \varphi_\theta:=\chi_\theta(\varphi).$$
When $|\theta|$ decreases to 0, $\chi_\theta$ decreases to
$\chi_0=\id$ and $\varphi_\theta$
decreases to $\varphi$. The following lemma gives some properties
of $\varphi_\theta$ where the 
coordinates $(x_1',y,v')$ are introduced in Section \ref{section_kahler}.

\begin{lemma} \label{lemma_varphi_theta}
There is a constant $c>0$ such that for $\theta\in\C^*$ small enough, 
$\ddc\varphi_\theta+\gamma$ vanishes on
$\{|x_1|>c|\theta|\}$. Moreover, we can write 
$$\ddc\varphi_\theta+\gamma=Adx_1\wedge d\overline x_1+B$$
where $A$ is a smooth function  such that $\|A\|_\infty\leq
c|\theta|^{-2}$ and $B$ is a smooth form such that $\|B\|_\infty\leq c|\theta|^{-1}$. 
\end{lemma}
\proof
Since $\widehat\Delta$ is
given by $x_1=0$ and $\ddc\varphi=[\widehat\Delta]-\gamma$, the function
$\psi:=\varphi-\log|x_1|$ is smooth.
By definition, $\varphi_\theta=\varphi$ on $\{\varphi>
\log|\theta|\}$ which contains $\{|x_1|>c|\theta|\}$ for some
constant $c>0$ large enough. So, we have for $|x_1|>c|\theta|$
$$\ddc\varphi_\theta+\gamma=\ddc\varphi+\gamma=0.$$
This proves the first assertion of the lemma.

For the second assertion, observe that $\varphi_\theta$ is constant on
$\{|x_1|< c'|\theta|\}$ for some constant $c'>0$. Therefore, it is enough to
consider the problem on the domain $\{c'|\theta|\leq |x_1|\leq
c|\theta|\}$ where we have $\ddc\varphi=\gamma$.
Observe that $\|\varphi\|_{\Cc^1}\lesssim |x_1|^{-1}$,
$\|\varphi\|_{\Cc^2}\lesssim |x_1|^{-2}$ and that the
derivatives of $\chi_\theta$ are bounded by a constant independent of
$\theta$. We have
$$\ddc\varphi_\theta=\ddc\chi_\theta(\varphi)=\chi_\theta''(\varphi)d\varphi\wedge\dc\varphi
+\chi_\theta'(\varphi)\ddc\varphi.$$
The last term is bounded. For the first term, we have since $\psi$ is smooth
$$d\varphi\wedge \dc\varphi=d(\log|x_1|+\psi)\wedge
\dc(\log|x_1|+\psi) = {i\over\pi}|x_1|^{-2}dx_1\wedge d\overline x_1 +
O(|x_1|^{-1}).$$
This implies the result. 
\endproof

\begin{lemma} \label{lemma_qpsh_theta_z}
The function $(\theta,z)\mapsto \varphi_\theta(z)$ can
  be extended to a
  quasi-p.s.h. function on $\C\times \widehat{X\times X}$ which is continuous
  outside $\{0\}\times\widehat\Delta$ and d.s.h. on
  $\P^1\times\widehat{X\times X}$. We have
$\varphi_0(z)=\varphi(z)$ and 
$\ddc \varphi_\theta(z)\geq -\lambda\widehat\omega(z)$
  on $\C\times \widehat{X\times X}$ for some constant
  $\lambda>0$. Moreover, $\ddc\varphi_\theta(z)$ can be written as a
  difference of positive closed currents on
  $\P^1\times\widehat{X\times X}$ which are smooth on
  $\C^*\times\widehat{X\times X}$. 
\end{lemma}
\proof
If $\psi(\theta,z):=\varphi(z)-\log|\theta|$, then we have on
$\C^*\times \widehat{X\times X}$
\begin{eqnarray*}
\ddc\varphi_\theta(z) & = & [\chi'(\psi)]^2d\psi\wedge
\dc \psi+\chi''(\psi)\ddc\psi\\
& \geq & \chi''(\psi)\ddc\psi=\chi''(\psi)\ddc\varphi(z).
\end{eqnarray*}
Hence, $\ddc\varphi_\theta(z)\geq -\lambda\widehat\omega(z)$, $\lambda>0$, on $\C^*\times \widehat{X\times X}$
because $\chi''$ is positive bounded and $\varphi$ is quasi-p.s.h. On
the other hand, by
definition, $\varphi_\theta(z)=\log|\theta|-1$ when $|\theta|\geq 1$
and $\varphi_\theta(z)$ is bounded from above when $|\theta|\leq 1$. 
By classical properties of p.s.h. functions,
$\varphi_\theta(z)$ can be extended to a quasi-p.s.h. function and the
estimate $\ddc\varphi_\theta(z)\geq -\lambda\widehat\omega(z)$ holds on 
$\C\times \widehat{X\times X}$. 

Since
$\varphi_\theta(z)=\log|\theta|-1$ for $|\theta|\geq 1$,
$\varphi_\theta(z)$ is d.s.h. on $\P^1\times\widehat{X\times X}$. We
have for the $\ddc$ operator on $\P^1\times\widehat{X\times X}$
$$\ddc \varphi_\theta(z)\geq -\big[\{\infty\}\times\widehat{X\times
  X}\big]-\lambda\widehat\omega(z).$$
So, we can write
$\ddc\varphi_\theta(z)$ as the following difference of two positive closed
currents
$$\Big(\ddc\varphi_\theta(z)+\big[\{\infty\}\times\widehat{X\times
  X}\big]+\lambda\widehat\omega(z)\Big) -\Big(\big[\{\infty\}\times\widehat{X\times
  X}\big]+\lambda\widehat\omega(z)\Big).$$
These currents are smooth on $\C^*\times \widehat{X\times X}$ since
$\varphi_\theta(z)$ is smooth there.

It remains to study $\varphi_\theta(z)$ when $\theta=0$ or
$\theta\rightarrow 0$. For $a\not\in\widehat\Delta$, we have
$\varphi_\theta(z)\rightarrow \varphi(a)$ when $(\theta,z)\rightarrow
(0,a)$. Therefore, $\varphi_\theta(z)$ is continuous out of
$\{0\}\times\widehat\Delta$ and $\varphi_0(z)=\varphi(z)$ outside
$\widehat\Delta$. Finally, since $\varphi_\theta\leq \max(\varphi,
\log|\theta|)$, we have that $\varphi_\theta(z)\rightarrow -\infty$
when $(\theta,z)$ tends to $\{0\}\times\widehat\Delta$. Since
$\varphi_\theta(z)$ is quasi-p.s.h. on $\C\times\widehat{X\times X}$, we deduce that
$\varphi_0(z)=-\infty=\varphi(z)$ 
on $\widehat\Delta$. 
\endproof

Define for $\theta\in\C^*$ the current $\Theta_\theta$ on
$\widehat{X\times X}$ by
$$\Theta_\theta:=(\ddc
\varphi_\theta+\gamma)\wedge \eta\quad
\mbox{and}\quad \Theta_0:=[\widehat\Delta]\wedge\eta,$$
see Example \ref{example_id}. 
Observe that $\Theta_\theta=\gamma\wedge \eta$ for $|\theta|>1$ since
in this case $\varphi_\theta$ is constant. So, define also
$\Theta_\infty:=\gamma\wedge \eta$. Since $\gamma$ and $\eta$ are
smooth, they can be written as differences of smooth positive closed
forms. By Lemma \ref{lemma_qpsh_theta_z}, for each $\theta\not=0$, $\ddc\varphi_\theta$ can be
written as a
difference of smooth positive closed forms on $X$ with masses bounded by a
constant independent of $\theta$. Therefore, 
we can write $\Theta_\theta:=\Theta_\theta^+-\Theta_\theta^-$ where
$\Theta_\theta^\pm$ are smooth positive closed with mass bounded by a
constant independent of $\theta$. We see that the family
$(\Theta_\theta)_{\theta\in\P^1}$ is continuous with respect to the
$\ast$-topology. Define
$$\L_\theta^\pm:=\L_{\Theta_\theta^\pm}\quad \mbox{and}\quad
\L_\theta:=\L_{\Theta_\theta}=\L_\theta^+-\L_\theta^-.$$
Note that $\L_\theta$ is symmetric.

\medskip

\noindent
$\bullet$ {\bf Regularization of currents.}
Regularization of positive closed currents on complex manifolds was developed by
Demailly in the case of bidegree $(1,1)$ \cite{Demailly3}. The case of
bidegree $(p,p)$ was studied in \cite{DinhSibony3}.
Define $S_\theta:=\L_\theta(S)$ for all currents $S$ in $\Dc_p$. 

\begin{lemma} \label{lemma_continuity_line}
The current $S_\theta$ depends
  continuously on $(\theta,S)$ for the $\ast$-topology on $S$,
  $S_\theta$. In particular, we have $\|S_\theta\|_*\leq c\|S\|_*$ for some constant
  $c>0$ independent of $S$ and $\theta$. 
Moreover, we have $\dist_2(S_\theta,S)\leq c|\theta|\|S\|_*$ with
$c>0$ independent of $S$ and $\theta$.
\end{lemma}
\proof
The estimate $\|S_\theta\|_*\leq c\|S\|_*$ 
is clear for $\theta=0$ since $\L_0=\id$,
see Example \ref{example_id}.
The case $\theta\not=0$ is a consequence of Lemma \ref{lemma_tr_norm}
applied to $\L_\theta^\pm$. 
When $\theta$ tends to $a\in\C^*$, then $\varphi_\theta$
converges in the $\Cc^\infty$-topology to $\varphi_a$. Therefore, $S_\theta$ depends
  continuously on $(\theta,S)$ for $\theta \in\C^*$.

It remains to prove the estimate on $\dist_2(S_\theta,S)$ for
$|\theta|\leq 1$. This and the triangle inequality
imply the continuity of $S_\theta$ at $\theta=0$, see also Proposition
\ref{prop_compare_dist}. We can assume that
$S$ is positive and that $\|S\|\leq 1$. 
Let $\Phi$ be a test form such
that $\|\Phi\|_{\Cc^2}\leq 1$. 
We have using the description of $\L_0$ in Example \ref{example_id}
\begin{eqnarray*}
\langle S_\theta-S,\Phi\rangle & = & \langle
\ddc(\varphi_\theta-\varphi)\wedge\eta\wedge
\Pi_1^*(S),\Pi_2^*(\Phi)\rangle\\
& = & \langle (\varphi_\theta-\varphi)\eta\wedge \Pi_1^*(S),
\Pi_2^*(\ddc\Phi)\rangle\\
&= & \big\langle S, (\Pi_1)_*\big((\varphi_\theta-\varphi)\eta\wedge\Pi_2^*(\ddc\Phi)\big)\big\rangle.
\end{eqnarray*}
We have to bound the last integral by $c|\theta|$ for some constant $c>0$.

Since $\|S\|\leq 1$, it is enough to show that the form
$(\Pi_1)_*\big[(\varphi_\theta-\varphi)\eta\wedge\Pi_2^*(\ddc\Phi)\big]$
has a $\|\cdot\|_\infty$-norm bounded by $c|\theta|$. The map $\Pi_1$
is a submersion. So, the coefficients of the considered form are equal
to some integrals of coefficients of
$(\varphi_\theta-\varphi)\eta\wedge\Pi_2^*(\ddc\Phi)$ on fibers of
$\Pi_1$. The  $\|\cdot\|_\infty$ estimate is not difficult to
obtain. Indeed, $\eta\wedge\Pi_2^*(\ddc\Phi)$ is a smooth form with
bounded $\|\cdot\|_\infty$-norm, the function $\varphi_\theta-\varphi$ has support in
a neighbourhood of $\widehat\Delta$ of size $\lesssim |\theta|$ and
satisfies  $|\varphi_\theta-\varphi|\lesssim
-\log\dist(\cdot,\widehat\Delta)$ near $\widehat\Delta$.
\endproof 

We deduce the following
result obtained in \cite{DinhSibony3}, 
see also Propositions \ref{prop_structural_line} and
\ref{prop_regularization_sp} below.

\begin{theorem} \label{th_regularization}
Smooth forms are dense in $\Dc_p$ and in $\Dc^0_p$ for the
$\ast$-topology. Moreover, there is a constant $c>0$ such that for
every current $S\in\Dc_p$, we can
write $S=S^+-S^-$ with $S^\pm\in\Cc_p$, $\|S^\pm\|\leq c\|S\|_*$ and
$S^\pm$ approximable by smooth forms in $\Cc_p$. 
\end{theorem} 
\proof
We prove the first assertion.
If $S$ is in $\Dc_p$, we can add to $S$ a smooth form in
order to obtain a current in $\Dc^0_p$. So, it is enough to approximate
currents $S$ in $\Dc^0_p$ by smooth forms in $\Dc^0_p$.
Observe that the problem is easy when $S$ is a form
of class $\Cc^1$. Indeed, we can write $S=\ddc U$ with $U$ of class
$\Cc^2$ and approximate $S$ uniformly by $S_\epsilon:=\ddc U_\epsilon$
where $U_\epsilon$ is smooth and $U_\epsilon\rightarrow U$ in the $\Cc^2$
topology.
It remains to approximate $S$ by $\Cc^1$ forms in $\Dc^0_p$.
Consider non-zero complex numbers $\theta_1,\ldots,\theta_{k+2}$. 
The currents $\Theta_{\theta_i}$ are smooth, then the associated
transforms $\L_{\theta_i}$ are regular.
By Lemma \ref{lemma_continuity_line}, we
can choose $\theta_i$ converging to 0 such that 
$\L_{\theta_{k+2}}\circ\cdots\circ \L_{\theta_1}(S)$
  converges to $S$. By Corollary \ref{cor_regular} and Lemma \ref{lemma_tr_norm},
  $\L_{\theta_{k+2}}\circ\cdots\circ \L_{\theta_1}(S)$ is a $\Cc^1$
    form in $\Dc_p^0$. This completes the proof of the first assertion.

For the second assertion, we can assume that $S$ is positive. Recall
that when $\theta\not=0$, $\L_\theta=\L_\theta^+-\L_\theta^-$ where
$\L_\theta^\pm$ are associated to smooth positive forms
$\Theta_\theta^\pm$ with mass bounded by a constant. Therefore,
$\L_{\theta_{k+2}}\circ\cdots\circ \L_{\theta_1}(S)$ is a difference
of $\Cc^1$ positive closed forms of bounded mass. We obtain the result
by extracting subsequences of forms converging to some currents
$S^\pm$.  
\endproof
 
\begin{corollary} \label{cor_mass_open}
Let $S$ be a current in $\Dc_p$ and $S'$ a current
  in $\Dc_{p'}$ with $p+p'\leq k$. Assume that $S$ restricted to an
  open set $W$ is a continuous form. Then $S\wedge S'$ is defined on
  $W$ and its mass on $W$ satisfies
$$\|S\wedge S'\|_W\leq c\|S\|_*\|S'\|_*$$
for some constant $c>0$ independent of $S$ and $S'$.
\end{corollary}
\proof
It is clear that $S\wedge S'$ is well-defined on $W$ and depends
continuously on $S'$ for the $\ast$-topology on $S'$. Therefore,
by Theorem
\ref{th_regularization}, we can assume that $S'$ is positive and
smooth. Now $S\wedge S'$ is defined on $X$ for every $S$ smooth or
not. We can assume that $S$ is positive but we may loose the
continuity of $S$. The current $S\wedge S'$ is
positive on
$X$. Its mass can be computed cohomologically. Therefore, we have
$$\|S\wedge S'\|_W\leq \|S\wedge S'\|\leq c\|S\|\|S'\|.$$
This implies the result.
\endproof

\begin{lemma} \label{lemma_deformation_c}
Let $S$ and $S_\theta$ be as above. Assume that $S$ is smooth.
Then $S_\theta$ is smooth for every $\theta$ and
$$\|S_\theta-S\|_\infty\leq c|\theta|\|S\|_{\Cc^1}$$
where $c>0$ is a constant independent of $S$ and $\theta$.
\end{lemma}
\proof
The current $S_0$ is equal to $S$. Hence, $S_0$ is smooth. For
$\theta\not=0$, $\L_\theta$ is a regular transform. Using the fact
that $\Pi_2$ is a submersion, we deduce easily that
$S_\theta=\L_\theta(S)$ is smooth. It remains to prove the estimate in
the lemma.

Assume that $\|S\|_{\Cc^1}\leq 1$. 
Observe that $(\Pi_2)_*(\Theta_\theta)$ is a closed current of
bidegree $(0,0)$ on $X$. So, it is defined by a constant function. On
the other hand, since $\Theta_\theta$ is cohomologous to $\Theta_0$,
$(\Pi_2)_*(\Theta_\theta)$ is cohomologous to
$$(\Pi_2)_*(\Theta_0)=(\pi_2)_*\pi_*(\Theta_0)=(\pi_2)_*[\Delta]=[X].$$
Hence, $(\Pi_2)_*(\Theta_\theta)$ is equal to 1. We deduce that 
$S=(\Pi_2)_*(\Pi_2^*(S)\wedge\Theta_\theta)$ and then 
$$S_\theta-S=(\Pi_2)_*(\pi^*(S')\wedge\Theta_\theta),\quad
\mbox{where}\quad S':=\pi_1^*(S)-\pi_2^*(S).$$
Observe that $\|S'\|_{\Cc^1}$ is bounded and the restriction
of $S'$ to $\Delta$ vanishes. If $S'':=\pi^*(S')$, then
$\|S''\|_{\Cc^1}$ is bounded and $S''$ restricted to $\widehat\Delta$
vanishes. Therefore, in the local coordinates near $\widehat\Delta$ as
in Section \ref{section_kahler}, we have
$$S''(x_1,y,v')=|x_1|A+Bdx_1+Cd\overline x_1$$
where $A$, $B$, $C$ are bounded forms.

The coefficients of $S_\theta-S$ at a point $y^0$ is computed by some
integrals involving the coefficients of
$S''\wedge\Theta_\theta=S''\wedge (\ddc\varphi_\theta+\gamma)\wedge\eta$ on
$\{y=y^0\}$. The above description of $S''$ together with Lemma
\ref{lemma_varphi_theta}, implies that these coefficients are
$\lesssim |\theta|$. The result follows.
\endproof

\noindent
$\bullet$ {\bf Green potential and $\ddc$-equation.}
Consider a current $S$ in $\Dc^0_p$ with $p\geq 1$. Then, since $[S]=0$, by
$\ddc$-lemma \cite{Demailly3, Voisin}, there is a real current $U_S$ of bidegree
$(p-1,p-1)$ such that $\ddc U_S=S$. We call $U_S$ {\it a potential} of
$S$. In order to construct an explicit potential and to estimate its norm,
we use a transform of bidegree $(-1,-1)$.
Choose a real smooth $(k-1,k-1)$-form
$\beta$ on $\widehat{X\times X}$ such that $\ddc\beta=-\gamma\wedge \eta+\pi^*(\alpha_\Delta)$ where
$\alpha_\Delta$, $\gamma$ and $\eta$ are introduced in Section
\ref{section_kahler}.
We can choose $\beta$ symmetric.
Consider the symmetric transform $\L_K$ with $K:=\varphi\eta-\beta$. The
following result was obtained in \cite[Proposition 2.1]{DinhSibony4},
see also \cite{GilletSoule}.

\begin{proposition} Let $S$ be a current in $\Dc^0_p$ with $p\geq
  1$. Then $U_S:=\L_K(S)$ is a potential of $S$. Moreover, we have
  $$\|U_S\|_{L^{1+1/k}}\leq c\|S\|_*$$ 
for some constant $c>0$ independent of  $S$. 
\end{proposition}
\proof
By Proposition \ref{prop_regular},
$S\mapsto\L_K(S)$ is continuous with respect to the 
$\ast$-topology on $S\in \Dc^0_p$ and the estimate on
$\|U_S\|_{L^{1+1/k}}$ is clear. We show that $\ddc U_S=S$. By Theorem
\ref{th_regularization}, 
it is enough to consider 
$S$ smooth. Define
$K':=\pi_*(K)$. This is a form smooth outside $\Delta$. 
We have seen in Section \ref{section_transform} that
$$|K'|\lesssim
-\log\dist(\cdot,\Delta)\dist(\cdot,\Delta)^{2-2k}$$
near $\Delta$. We also have 
$$\ddc K'=\pi_*(\ddc
K)=\pi_*\big([\widehat\Delta]\wedge\eta-\pi^*(\alpha_\Delta)\big)
=[\Delta]-\alpha_\Delta.$$ 
So, $K'$ is
a kernel for solving the $\ddc$-equation on $X$. Since $S$ is smooth, we have
$U_S=(\pi_2)_*(\pi_1^*(S)\wedge K')$ and
\begin{eqnarray*}
\ddc U_S & = & (\pi_2)_*\big(\pi_1^*(S)\wedge [\Delta]\big) -
(\pi_2)_*\big(\pi_1^*(S)\wedge \alpha_\Delta\big) \\
& = & S- (\pi_2)_*\big(\pi_1^*(S)\wedge
\alpha_\Delta\big)=S,
\end{eqnarray*}
where the last identity is obtained using that $[S]=0$ and 
that $\alpha_\Delta$ is a combination of forms of type
$\beta(x)\wedge\beta'(y)$ with $\beta$ and $\beta'$ closed.
\endproof

\begin{definition}\rm
We call $\L_K(S)$ {\it the Green potential of $S$}. 
\end{definition}

Note that the Green potential depends on the choice of $K$.

\begin{remark}\rm
The transform associated to ${i\over \pi}\partial K$ solves the $\dbar$-equation
on $X$ and ${i\over\pi}\pi_*(\partial K)$ is a kernel which solves the $\dbar$-equation.
\end{remark}


\section{Structural varieties and super-potentials} \label{section_sp}

In this section, we define for each positive closed $(p,p)$-current a
super-potential which is a function on the space of smooth forms in 
$\Dc^0_{k-p+1}$. In this space, we construct some special structural lines
parametrized by the projective line $\P^1$. The restriction of the
super-potential to such a structural line is a d.s.h. function. This
is a key point in our study. We will also consider currents with regular
super-potentials and their intersection.

\subsection{Structural varieties in the space of currents} \label{section_structural}

Consider 
a current $\Rc$ on $V\times
X$ which is a difference of two positive closed $(s,s)$-currents. We
will use in next sections the case where $s=k-p+1$. By
Proposition \ref{prop_slice}, for $\theta$ in $V$ outside a locally pluripolar set,
the slice $R_\theta:=\langle \Rc,\pi_V,\theta\rangle$ is well-defined and is a
current in $\Dc_s$. Its cohomology class does not depend on $\theta$. 
Assume that $R_\theta$ is in $\Dc^0_s$,
i.e. $[R_\theta]=0$.
So, we obtain
a map $\tau:V\rightarrow\Dc^0_s$ given by $\theta\mapsto R_\theta$ which
is defined out of a locally pluripolar set. 

\begin{definition}\rm
We say that $\tau$ or the family
$(R_\theta)_{\theta\in V}$ defines a {\it structural variety} of $\Dc^0_s$. When
$R_\theta$ is defined for every $\theta$ and depends continuously on $\theta$ for the $\ast$-topology, we say that the
structural variety is {\it continuous}. 
\end{definition}

In what follows, we use some {\it structural lines},
i.e. structural varieties parametrized by the projective line $\P^1=\C\cup\{\infty\}$. 
Let $\L_\theta:=\L_{\Theta_\theta}$ be transforms defined in
Section \ref{section_regularization}.
For a given current $R$ in $\Dc^0_s$ and for $\theta\in \C\cup\{\infty\}$, consider the
current $R_\theta:=\L_\theta(R)$.
Recall that $\L_\theta$, $\Theta_\theta$, $R_\theta$ do not depend on
$\theta$ when $|\theta|\geq 1$ and that $\L_0=\id$, $R_0=R$.

\begin{proposition} \label{prop_structural_line}
The family of currents $(R_\theta)_{\theta\in\P^1}$ defines a
  continuous structural line in $\Dc^0_s$ which
  depends linearly on $R$. Moreover,  there is a constant $c>0$
  independent of $R$ such that $\|R_\theta\|_*\leq c\|R\|_*$ for every
  $\theta$. 
\end{proposition}

\begin{definition}\rm
We call  $(R_\theta)_{\theta\in\P^1}$ {\it the special structural
  line} associated to $R$.
\end{definition}

\noindent
{\bf Proof of Proposition \ref{prop_structural_line}.}
The linear dependence on $R$ is clear. The continuity of
$(R_\theta)_{\theta\in\P^1}$ and the
estimate on
 $\|R_\theta\|_*$ are proved in Lemma \ref{lemma_continuity_line}.
Let $\tau_0$ denote the projection of $\P^1\times \widehat{X\times X}$
on $\P^1$ and $\tau$ the projection on $\widehat{X\times X}$. Consider
$\varphi_\theta(z)$ as a function on $\C\times \widehat{X\times
  X}$. Define a current $\widehat\Rc$ on $\C^*\times\widehat{X\times
  X}$ by 
$$\widehat\Rc(\theta,z)  :=  \big(\ddc\varphi_\theta(z)+\tau^*(\gamma)\big)\wedge
\tau^*(\eta)\wedge \tau^*(\Pi_1^*(R)).$$
In this wedge-product, each current is a difference of positive closed
currents with bounded mass in $\P^1\times \widehat{X\times X}$. We can
apply Corollary \ref{cor_mass_open} to the current $\widehat\Rc$, which is
well-defined on $\C^*\times \widehat{X\times X}$, and Skoda's theorem
\cite{Skoda} on the extension of positive closed currents. Hence,
the trivial
extension of $\widehat\Rc$ is a difference of positive closed currents on
$\P^1\times\widehat{X\times X}$ with bounded mass. Denote also by $\widehat \Rc$ this extension.

On $\C^*\times\widehat{X\times X}$, in the definition of $\widehat \Rc$, all currents
except $R$, are smooth. We deduce easily from the slicing theory that
$$\langle\widehat\Rc,\tau_0,\theta\rangle=(\ddc_z\varphi_\theta+\gamma)\wedge\eta\wedge
\Pi_1^*(R)$$
where we identify $\{\theta\}\times\widehat{X\times X}$ with
$\widehat{X\times X}$. Let $\tau_2:=(\tau_0,\Pi_2\circ\tau)$ denote
the projection of $\P^1\times\widehat{X\times X}$ onto the product of $\P^1$
with the second factor $X$. Define $\Rc:=(\tau_2)_*(\widehat\Rc)$. It is
deduced from the slicing theory that
$$\langle
\Rc,\pi_{\P^1},\theta\rangle=(\Pi_2)_*\langle\widehat\Rc,\tau_0,\theta\rangle
= R_\theta,$$
for $\theta\in\C^*$. Recall that $R_\theta$ depends continuously on
$\theta\in\P^1$. By 
Remark \ref{rk_slice_continuity}, the identity $\langle
\Rc,\pi_{\P^1},\theta\rangle= R_\theta$
holds for any
$\theta\in\P^1$. So, $(R_\theta)_{\theta\in\P^1}$ is a structural line.
\hfill $\square$

\begin{remark}\rm
We can prove that $\theta\mapsto \L_\theta^{k+2}(R)$ defines a
continuous 
structural line. In this case, for $\theta\not=0$,
$\L_\theta^{k+2}(R)$ is a $\Cc^1$ form.
\end{remark}

\subsection{Super-potentials of currents} \label{section_sp_def}

Consider a current $S$ in $\Dc_p$. The
super-potentials of $S$ are defined (at least) on the smooth forms in
$\Dc^0_{k-p+1}$. They are unique under 
apropriate normalization.

Let $\alpha=\{\alpha_1,\ldots,\alpha_h\}$ with $h:=\dim H^{p,p}(X,\R)$
be a fixed family of real smooth
closed $(p,p)$-forms such that the family of classes
$[\alpha]=\{[\alpha_1],\ldots,[\alpha_h]\}$ is a basis of
$H^{p,p}(X,\R)$.  We can find  a family
$\alpha^\vee=\{\alpha_1^\vee,\ldots,\alpha_h^\vee\}$ of real smooth
closed $(k-p,k-p)$-forms such that
$[\alpha^\vee]=\{[\alpha_1^\vee],\ldots,[\alpha_h^\vee]\}$ is the dual
basis of $[\alpha]$ with respect to the cup-product $\smile$. 
Let $R$ be a current in $\Dc^0_{k-p+1}$ and $U_R$ a potential of
$R$. Adding to $U_R$ a suitable combination of $\alpha_i^\vee$ allows to assume
that $\langle U_R,\alpha_i\rangle=0$ for
$i=1,\ldots,h$. We say that $U_R$ is {\it $\alpha$-normalized}. 

\begin{lemma} \label{lemma_integral_sp}
Assume that $S$ is smooth or that $R$, $U_R$ are smooth. Then $\langle
S,U_R\rangle$ does not depend on the choice of $U_R$.  
Assume that $[S]=0$. Let $U_S$ be a potential of $S$, smooth if $S$ is
smooth. Let $U_R'$ be another potential of $R$, smooth when $R$ is
smooth. Then $\langle S,U_R'\rangle=\langle S,U_R\rangle=\langle
U_S,R\rangle$. In particular, 
$\langle S,U_R\rangle$ does not depend on $\alpha$ and $\alpha^\vee$.
\end{lemma}
\proof
Let $U_R'$ be another $\alpha$-normalized potential of $R$. We have
$\ddc (U_R-U'_R)=0$ and $[\alpha_i]\smile [U_R-U_R']=0$ for every $i$. Since
$[\alpha]$ is a basis of $H^{p,p}(X,\R)$, we deduce that $[S]\smile
[U_R-U_R']=0$. Hence, $\langle S,U_R\rangle=\langle
S,U_R'\rangle$. So, $\langle S, U_R\rangle$ does not depend on the
choice of $U_R$. If $[S]=0$, we have 
$$\langle S,U_R'\rangle =\langle \ddc U_S,U_R'\rangle =\langle U_S,\ddc
U_R'\rangle = \langle U_S,R\rangle.$$
These identities hold for all $U_R'$ not necessarily normalized, in particular for $U_R$.
Note that the smoothness of $S$, $U_S$ or $R$, $U_R$, $U_R'$ implies that
the considered integrals are meaningful.
\endproof

\begin{definition}\rm
{\it The $\alpha$-normalized super-potential $\Uc_S$ of $S$}, is
the following function defined on smooth forms
$R$ in $\Dc^0_{k-p+1}$ by 
\begin{equation} \label{eq_sp_def}
\Uc_S(R):=\langle S,U_R\rangle,
\end{equation}
where $U_R$ is an $\alpha$-normalized smooth potential of $R$. We say
that $S$ has {\it a continuous super-potential}\footnote{this 
  is equivalent to the notion of PC current introduced in \cite{DinhSibony4}.} if $\Uc_S$ can be
extended to a function on $\Dc^0_{k-p+1}$ which is continuous with respect to the
$\ast$-topology. In this case, the extension is also denoted by $\Uc_S$
and is also called super-potential of $S$.
\end{definition}

Note that the $\alpha$-normalized super-potential of $\alpha_i$ is
identically zero. By Lemma \ref{lemma_integral_sp}, when $[S]=0$, the super-potential $\Uc_S$
does not depend on the choice of $\alpha$. When $S$ is smooth then $S$
has a continuous super-potential and the formula (\ref{eq_sp_def})
holds for all $R$ in $\Dc_{k-p+1}^0$. In this case, if $[S]=0$ and if
$U_S$ is a smooth potential of $S$, we also have $\Uc_S(R)=\langle
U_S,R\rangle$. 

\begin{proposition} \label{prop_current_sp}
Let $S$ and $S'$ be two currents in $\Dc_p$ such that
$[S]=[S']$. 
If they have the same $\alpha$-normalized super-potential then they
are equal.
\end{proposition}
\proof
The $\alpha$-normalized super-potential $\Uc_{S''}$ of $S'':=S-S'$
vanishes identically. 
If $U$ is a real smooth $(k-p,k-p)$-form, then $U$ is a potential of
$\ddc U$ which is a form in $\Dc^0_{k-p+1}$. Since $[S'']=0$, it follows
from Lemma \ref{lemma_integral_sp}, that $\langle
S'',U\rangle=\Uc_{S''}(\ddc U)=0$. Hence, $S''=0$.
\endproof

Here is one of the fundamental properties of super-potential. It can be
extended to more general structural varieties but, for simplicity we restrict ourselves
to this particular case.

\begin{proposition} \label{prop_sp_line}
Let $(R_\theta)_{\theta\in\P^1}$ be the special structural line associated
  to a smooth form $R\in\Dc^0_{k-p+1}$. Let $S$ be a current in
  $\Dc_p$. 
Then $\theta\mapsto \Uc_S(R_\theta)$ is a continuous d.s.h. function on
$\P^1$ which is constant on $\{|\theta|\geq 1\}$. Moreover, we have
$$\|\ddc_\theta\Uc_S(R_\theta)\|_*\leq c\|S\|_*\|R\|_*$$
where $c>0$ is a constant independent of $R$
and $S$.
\end{proposition}
\proof
By Lemma \ref{lemma_deformation_c} applied to $R_\theta$, the function 
$H(\theta):=\Uc_S(R_\theta)$ is continuous on
$\P^1$. It remains to bound the mass of
$\ddc H$. Since this function depends continuously
on $S$, by Theorem \ref{th_regularization}, we can assume that $S$ is
smooth.
Recall that the $\alpha$-normalized super-potential of $\alpha_i$ is zero.
Subtracting from $S$ a combination of $\alpha_i$ allows to assume
that $[S]=0$. So, we can use the last assertion of Lemma
\ref{lemma_integral_sp}: if $U$ is a smooth potential of $S$, then
$H(\theta)=\langle U,R_\theta\rangle$. 

It is enough to estimate the mass of $\ddc H$ on $\C^*$. Indeed, the
continuity of $H$ implies that $\ddc H$ has no mass on  finite
sets. Consider in $\P^1\times\widehat{X\times X}$ the currents
$$\widehat \Rc_U:=\widehat \Rc\wedge
\tau^*\Pi_2^*(U)\quad\mbox{and}\quad \widehat \Rc_S:=\widehat \Rc\wedge
\tau^*\Pi_2^*(S).$$
These currents are smooth on $\C^*\times\widehat{X\times X}$. 
A direct computation gives $H=(\tau_0)_*(\widehat\Rc_U)$ and $\ddc H=(\tau_0)_*(\widehat
\Rc_S)$ on $\C^*$. So, it is enough to estimate the mass of $\widehat \Rc_S$ on
$\C^*\times\widehat{X\times X}$. By Corollary \ref{cor_mass_open},
since $S$ and $R$ are smooth, this mass is bounded by a constant times
$$\|\widehat\Rc\|_*\|\tau^*\Pi_2^*(S)\|_*
\lesssim \|\tau^*\Pi_1^*(R)\|_*\|\tau^*\Pi_2^*(S)\|_*\lesssim \|R\|_*\|S\|_*,$$
where the last inequality follows from the fact that $\tau$, $\Pi_1$, $\Pi_2$ are submersions.
\endproof

\begin{lemma} \label{lemma_sp_theta}
Let $\Uc_{S_\theta}$ be the $\alpha$-normalized
  super-potential of $S_\theta$. If $[S]=0$, then $\Uc_{S_\theta}(R)=\Uc_S(R_\theta)$ for
  $R$ smooth. 
\end{lemma}
\proof
Since $S_\theta=\L_\theta(S)$, the Green potential of $S_\theta$ is
equal to $\L_K\L_\theta(S)$. Using the symmetry of $\L_\theta$ and
$\L_K$, we have by Lemma \ref{lemma_integral_sp}
$$\Uc_{S_\theta}(R)=\langle \L_K\L_\theta(S),R\rangle = \langle S,
\L_\theta\L_K(R)\rangle.$$
On the other hand, 
$$\ddc \L_\theta\L_K(R)=\L_\theta(\ddc\L_K(R))=\L_\theta(R)=R_\theta.$$
It follows that $\Uc_{S_\theta}(R)=\Uc_S(R_\theta)$.
\endproof

The following result is an analogue of the estimate in Proposition
\ref{prop_exp} for super-potentials of currents.

\begin{theorem} \label{th_main_estimate}
Let $S$ be a current in
  $\Dc_p$ and $\Uc_S$ the $\alpha$-normalized super-potential of
  $S$. Then we have for $R$ smooth in $\Dc^0_{k-p+1}$ with
  $\|R\|_*\leq 1$ 
$$|\Uc_S(R)|\leq c\|S\|_*\big(1+\log^+\|R\|_{\Cc^1}\big),$$
where  $\log^+:=\max(\log,0)$ and $c>0$ is a constant independent of $S$, $R$.
\end{theorem}
\proof
Subtracting from $S$ a combination of $\alpha_i$ allows to assume that
$[S]=0$. We can also assume that $\|S\|_*=1$. Let $U_S$ be the Green potential of $S$. 
By Lemma \ref{lemma_integral_sp}, $\Uc_S(R)=\langle
U_S,R\rangle$. 
We have to show that
$$M_{S,R}:={|\langle U_S,R\rangle|\over  1+\log^+\|R\|_{\Cc^1}}$$
is bounded when $\|R\|_*\leq 1$. The proof uses Proposition
\ref{prop_regular}, Lemma \ref{lemma_dsh_exp} and special structural lines in
$\Dc^0_{k-p+1}$. 
Consider the  numbers $q_n\geq 1$ given  by the induction identity
$q_n^{-1}=q_{n-1}^{-1}-1+(1+1/k)^{-1}$ for $n\leq k+1$ with
$q_0=1$. We have
$q_{k+1}=\infty$. 

\medskip

\noindent
{\bf Claim.} {\it For every $0\leq n\leq  k+1$ and $M>0$,  
there is a constant $c>0$
independent of $S$, $R$ such that $M_{S,R}\leq c$ if $\|R\|_*\leq 1$
and $\|R\|_{L^{q_n}}\leq M$.} 

\medskip

For $n=0$, the claim implies the theorem, i.e. $M_{S,R}$ is bounded when
$\|R\|_*\leq 1$. Indeed, we have $\|R\|_{L^1}\lesssim\|R\|\lesssim \|R\|_*$
and then the hypothesis  $\|R\|_{L^{q_0}}\leq M$ is satisfied.
We prove now the claim using a decreasing induction on $n$. 
For $n=k+1$, by Proposition \ref{prop_regular}, $\|U_S\|_{L^1}$ is
bounded uniformly on $S$. If $\|R\|_\infty$ is bounded, it is clear
that $\langle U_S,R\rangle$ is bounded. So, the claim is true for
$n=k+1$. Assume now that the claim is true for $n+1$. 
We check it for $n$ and we only have to
consider the case where $\|R\|_{\Cc^1}$ is large. 

Let $R_\theta$ be as
above and define $H_{S,R}(\theta):=\Uc_S(R_\theta)$. By Proposition
\ref{prop_sp_line}, $H_{S,R}$ is a continuous d.s.h. function on
  $\P^1$. It is equal to some constant $c_{S,R}$ 
on $\{|\theta|\geq 1\}$. We have $c_{S,R}=\langle
  U_S,R_\infty\rangle$. Moreover, $\|\ddc
  H_{S,R}\|$ is bounded uniformly on $S$, $R$. On the other hand, by Propositions \ref{prop_regular}
  and \ref{prop_structural_line}, $R_\infty$ is a smooth form in $\Dc^0_{k-p+1}$ with
  bounded $L^{q_{n+1}}$-norm and bounded $\|\cdot\|_*$-norms. Since 
  $\langle U_S,R\rangle$ depends linearly on $R$, we can apply the claim to $R_\infty$ and 
  deduce that
$$c_{S,R}\lesssim 1+\log^+\|R_\infty\|_{\Cc^1}\lesssim 1+\log^+\|R\|_{\Cc^1}.$$
Because $\langle U_S,R\rangle=H_{S,R}(0)$, 
it is enough to show that 
$$|H_{S,R}(0)-c_{S,R}|\lesssim 1+\log^+\|R\|_{\Cc^1}.$$ 
Since $\|\ddc H_{S,R}\|$ is uniformly bounded, by Lemma
\ref{lemma_dsh_exp} applied to $H_{S,R}-c_{S,R}$, there is a $\theta$
with $|\theta|\leq \|R\|_{\Cc^1}^{-1}$ such that 
$$|H_{S,R}(\theta)-c_{S,R}|\lesssim 1+\log^+\|R\|_{\Cc^1}.$$
On the other hand, Lemma \ref{lemma_deformation_c} implies that
$$\|R-R_\theta\|_\infty\lesssim |\theta|\|R\|_{\Cc^1}\leq 1.$$
Therefore, using that $\|U_S\|_{L^1}$ is bounded, we obtain
\begin{eqnarray*}
|H_{S,R}(0)-c_{S,R}| & \leq &  |H_{S,R}(0)-H_{S,R}(\theta)|
+|H_{S,R}(\theta)-c_{S,R}| \\
& = &  |\langle U_S,R-R_\theta\rangle|
+|H_{S,R}(\theta)-c_{S,R}| \\
& \lesssim &  1+\log^+\|R\|_{\Cc^1}.
\end{eqnarray*}
This completes the proof.  
\endproof

We will use the following notion of convergence.

\begin{definition} \rm
Let $(S_n)$ be a sequence of currents converging in $\Dc_p$ to a
current $S$. Let $\Uc_S$, $\Uc_{S_n}$ be the $\alpha$-normalized
super-potentials of $S$, $S_n$.
We say that the convergence is {\it SP-uniform} if 
$\Uc_{S_n}$ converge to $\Uc_S$
uniformly on any $\ast$-bounded set of smooth forms in $\Dc^0_{k-p+1}$.
\end{definition}

By linearity, it is enough to check the SP-uniform convergence on the
unit ball of $\Dc^0_{k-p+1}$.  
This notion does not depend on
$\alpha$. Indeed, by Lemma \ref{lemma_integral_sp}, the case 
where $[S_n]=[S]=0$ is clear. Since $[S_n]$ converge to $[S]$, we obtain
the general case by adding to
$S_n$ and $S$ suitable combinations of $\alpha_i$.
Moreover, if $S_n$ and $S$ have continuous
super-potentials, then since smooth forms are dense in $\Dc^0_{k-p+1}$,
the extensions of $\Uc_{S_n}$ converge to the extension of $\Uc_S$
uniformly on $\ast$-bounded subsets of $\Dc^0_{k-p+1}$.

\begin{proposition} \label{prop_regularization_sp}
Let $S$ be a current in $\Dc_p$ with continuous super-potentials. Then
$S_\theta$ has continuous super-potentials and $S_\theta$ converges SP-uniformly to
$S$ when $\theta\rightarrow 0$. In particular, $S$ can be approximated
SP-uniformly by smooth forms in $\Dc_p$.  
\end{proposition}
\proof
Observe that the second assertion is deduced from the first one as in the
proof of Theorem \ref{th_regularization}. We prove now the first assertion.
When $S$ is smooth, by Proposition \ref{prop_regular}, $S_\theta$ converges to $S$
in the $\Cc^1$-topology and the result is clear. Adding to $S$ a
combination of $\alpha_i$ allows to assume that $[S]=0$. Then, we also
have $[S_\theta]=0$ for every $\theta$. We only have to consider
$|\theta|\leq 1$ since $S_\theta$ does not depend on $\theta$ when
$|\theta|\geq 1$. 
Let $R$ be a current in $\Dc_{k-p+1}^0$ with $\|R\|_*\leq 1$. Let
$\Uc_{S_\theta}$ denote the super-potential of $S_\theta$.
By Lemma \ref{lemma_sp_theta}, when $R$ is smooth, we have
$\Uc_{S_\theta}(R)=\Uc_S(R_\theta)$. Since $R_\theta$ depends
continuously on $R$, $\Uc_{S_\theta}$ admits a continuous extension to
$\Dc_{k-p+1}^0$ and the last identity holds for all $R$. 

It remains to check that $\Uc_{S_\theta}$ converges SP-uniformly to
$\Uc_S$. Recall that since $\Uc_S$ is continuous, if $\|R\|_*$ is
bounded, we have $\Uc_S(R)\rightarrow 0$ when $\|R\|_{\Cc^{-2}}\rightarrow 0$.
We also have $\Uc_{S_\theta}(R)-\Uc_S(R)=\Uc_S(R_\theta-R)$. When
$\theta\rightarrow 0$ and $\|R\|_*\leq 1$, $\|R_\theta-R\|_*$ is bounded and by Lemma
\ref{lemma_continuity_line}, $\|R_\theta-R\|_{\Cc^{-2}}$ tends to 0 uniformly on $R$. 
Therefore, $\Uc_S(R_\theta-R)$ tends to 0 uniformly on $R$ with $\|R\|_*\leq 1$. The result follows.
\endproof

\subsection{Intersection of currents} \label{section_intersection}

We will define the intersection of two currents such that at least one
of them has a continuous super-potential. The theory of intersection
is far from being complete but we will see that the following
properties suffice in order to study the dynamics of automorphisms. We
refer the reader to \cite{Demailly2,Demailly1, FornaessSibony2} for the
theory of intersection with currents of bidegree $(1,1)$ and
\cite{DinhSibony7, DinhSibony10}
for the case of bidegree $(p,p)$ on local setting or on homogeneous
manifolds, see also \cite{BostGilletSoule,GilletSoule}.

Let $S$ be a current in $\Dc_p$ and $S'$ a current in $\Dc_{p'}$ with
$p+p'\leq k$. 
Assume that $S$ has a continuous super-potential. We will define the
intersection $S\wedge S'$ as a current in $\Dc_{p+p'}$. This
wedge-product satisfies some continuity properties. Let
$\Uc_S$ be the $\alpha$-normalized super-potential of $S$ and 
let $(a_1,\ldots,a_h)$ denote the coordinates of $[S]$ in the basis
$[\alpha]$. Define for any
test smooth real form $\Phi$ of bidegree  $(k-p-p',k-p-p')$:
$$\langle S\wedge S',\Phi\rangle:=\Uc_S(\ddc\Phi\wedge S') +
\sum_{1\leq i\leq h} a_i\langle \alpha_i, \Phi\wedge S'\rangle.$$

\begin{lemma} \label{lemma_wedge_smooth}
Assume that $S$ or $S'$ is smooth. Then $S\wedge S'$ coincides with the usual
wedge-product of $S$ and $S'$.
\end{lemma}
\proof
Assume that $S'$ is smooth. 
Observe that $\Phi\wedge S'$ is a potential of
$\ddc\Phi\wedge S'$. Define $m_i:=\langle\alpha_i,\Phi\wedge
S'\rangle$. Then $\Phi\wedge S'-\sum m_i\alpha_i^\vee$ is an
$\alpha$-normalized potential of $\ddc\Phi\wedge S'$. Therefore,
\begin{eqnarray*}
\Uc_S(\ddc\Phi\wedge S')+\sum a_im_i & = & \langle S,\Phi\wedge S'\rangle -\sum
m_i\langle S,\alpha_i^\vee\rangle+\sum a_im_i \\
& = & \langle S,\Phi\wedge S'\rangle.
\end{eqnarray*} 
This implies that $S\wedge S'$ coincides with the usual
wedge-product of $S$ and $S'$. The computation still holds when $S$ is
smooth but $S'$ is singular.
\endproof

\begin{theorem} \label{th_wedge}
Let $S$ be a current in $\Dc_p$ and $S'$ a current in $\Dc_{p'}$ with
$p+p'\leq k$. Assume that $S$ has continuous
super-potentials. Then $S\wedge S'$, defined as above, is a current in
$\Dc_{p+p'}$ which depends linearly on $S$, $S'$. Moreover, we have
$$[S\wedge S']=[S]\smile[S']\quad \mbox{and}\quad \|S\wedge S'\|_*\leq c\|S\|_*\|S'\|_*$$
for some constant $c>0$ independent of $S$, $S'$. Let
$S_n$ be  currents in $\Dc_p$ with continuous super-potentials
converging SP-uniformly to $S$ and $S_n'$ be  currents 
converging in $\Dc_{p'}$ to $S'$. Then, $S_n\wedge S'_n$ converge in
$\Dc_{p+p'}$ to $S\wedge S'$.
\end{theorem}
\proof
It is clear that $\langle S\wedge S',\Phi\rangle$ depends continuously
on the smooth test form $\Phi$. Hence, $S\wedge S'$ is a
current. Clearly, this current depends linearly on $S$ and $S'$. By
definition, since $\Uc_S$ is continuous, $S\wedge S'$ depends
continuously on $S'$. We deduce using Theorem \ref{th_regularization} that $[S\wedge S']=[S]\smile [S']$
since this identity holds for $S'$ smooth.
In order to estimate $\|S\wedge S'\|_*$, it is
enough to assume that $S'$ is smooth positive. Writing $S$ as a
difference of positive closed current, we see that $\|S\wedge
S'\|_*\lesssim \|S\|_*\|S'\|_*$. We use here that the mass of a positive closed
current depends only on its cohomology class. The last assertion of
the theorem is deduced directly from the definition of $S\wedge S'$. 
\endproof

\begin{proposition} \label{prop_cv_sp_wedge}
Assume that $S$, $S'$, $S_n$, $S_n'$ have
  continuous super-potentials and that $S_n$,
  $S_n'$ converge SP-uniformly to $S$, $S'$ respectively.
  Then $S\wedge S'$ and $S_n\wedge S_n'$ have also continuous super-potentials and 
$S_n\wedge S_n'$ converge SP-uniformly to $S\wedge S'$.
\end{proposition}
\proof
The proposition is clear when $S$, $S_n$ are linear combinations of
$\alpha_i$.
Subtracting from $S$ and $S_n$ suitable combinations of $\alpha_i$
allows to assume that $[S]=[S_n]=0$. 
So, if $\Phi$ is a smooth test form we have by definition
$\langle S\wedge S',\Phi\rangle = \Uc_S(S'\wedge \ddc\Phi)$. We deduce
that if $R$ is smooth then $\Uc_{S\wedge S'}(R)=\Uc_S(S'\wedge
R)$. Since $S$ and $S'$ have continuous super-potentials,
$\Uc_S(S'\wedge R)$ can be extended continuously to $R$ in
$\Dc_{k-p-p'+1}^0$. So, $S\wedge S'$ has a continuous
super-potential and the identity $\Uc_{S\wedge S'}(R)=\Uc_S(S'\wedge
R)$ holds for all $R$ in $\Dc_{k-p-p'+1}^0$. In the same way, we prove that $S_n\wedge S_n'$ has a
continuous super-potential and 
$\Uc_{S_n\wedge S'_n}(R)=\Uc_{S_n}(S'_n\wedge
R)$. It is now clear that $S_n\wedge S_n'$ converge SP-uniformly to
$S\wedge S'$.
\endproof

The following result shows that the wedge-product is commutative and associative.
The first property allows to define 
$S'\wedge S:=S\wedge S'$ when $S$ has continuous super-potentials and
$S'$ is singular. 

\begin{proposition} Let $S_i$, $i=1,2,3$, be currents in $\Dc_{p_i}$. 
Assume that $S_1$ and $S_2$ have continuous super-potentials,
  then 
$$S_1\wedge S_2=S_2\wedge S_1\quad \mbox{and}\quad (S_1\wedge S_2)\wedge S_3=S_1\wedge (S_2\wedge S_3).$$
\end{proposition}
\proof
The proposition is clear when $S_1$ and $S_2$ are smooth. The general
case is deduced from this particular case using Propositions
\ref{prop_regularization_sp}, \ref{prop_cv_sp_wedge} and Theorem \ref{th_wedge}.
\endproof

\begin{remark}\rm
Assume that $S$ and $S'$ are positive currents. By Theorem
\ref{th_wedge}, if $S$ is SP-uniformly 
approximable by smooth positive closed $(p,p)$-forms, then $S\wedge
S'$ is positive. This is also the case when $S'$ can be approximated
by positive smooth forms. In general, we don't know if $S\wedge S'$ is
always positive when $S$ or $S'$ has continuous super-potentials.
\end{remark}

\subsection{H{\"o}lder super-potentials and moderate currents}

Consider a current $S$ in $\Dc_p$ with continuous super-potentials. 
Its super-potentials are defined on $\Dc^0_{k-p+1}$. 

\begin{definition}\rm
We say that $S$ has  {\it a H{\"o}lder continuous super-potential} if it
admits a super-potential which is H{\"o}lder continuous on $\ast$-bounded subsets
of $\Dc^0_{k-p+1}$ with respect to $\dist_l$ for some real number $l>0$.
\end{definition}

In order to prove that $\Uc_S$ is H{\"o}lder continuous, it is enough to
show that $|\Uc_S(R)|\lesssim \|R\|_{\Cc^{-l}}^\lambda$ for
$\|R\|_*\leq 1$ and for some constant $\lambda>0$. 
By Proposition \ref{prop_compare_dist}, the definition does not depend on the choice of
$l$. One checks easily that the super-potentials of smooth forms are
H{\"o}lder continuous. Hence, if $S$ admits a H{\"o}lder continuous
super-potential, all the super-potentials of $S$ are H{\"o}lder
continuous. In other words, this notion does not depend on the
normalization of the super-potential.

\begin{proposition} \label{prop_holder_wedge}
Let $S$, $S'$ be currents in $\Dc_p$ and $\Dc_{p'}$, $p+p'\leq k$,
having H{\"o}lder continuous super-potentials. Then $S\wedge S'$ has a
H{\"o}lder continuous super-potential.
\end{proposition}
\proof
We can assume that $[S]=0$ and $[S']=0$. Let $\Uc_S$,
$\Uc_{S'}$ and $\Uc$ be the
super-potentials of $S$, $S'$ and $S\wedge S'$ respectively. So,
for $R$ in $\Dc^0_{k-p-p'+1}$ we have 
$\Uc(R)=\Uc_S(S'\wedge R)$. 
It is enough to prove for $R$ in a $\ast$-bounded subset of $\Dc^0_{k-p-p'+1}$ that
$$\|S'\wedge R\|_{\Cc^{-4}}\lesssim \|R\|_{\Cc^{-2}}^\lambda,$$
where $\lambda>0$ is a constant. Using a regularization, we can assume that $R$ is
smooth.
Let $U'$ be a potential of $S'$ and $\Phi$ a test form with
$\|\Phi\|_{\Cc^4}\leq 1$. We have since $\Uc_{S'}$ is H{\"o}lder
continuous and $\|\ddc\Phi\|_{\Cc^2}$ is bounded
\begin{eqnarray*}
\|S'\wedge R\|_{\Cc^{-4}} & = & \sup_\Phi |\langle S'\wedge R, \Phi\rangle|
=\sup_\Phi |\langle U'\wedge R, \ddc \Phi\rangle| \\
& = & \sup_\Phi |\langle U', R \wedge \ddc \Phi\rangle|
= \sup_\Phi |\Uc_{S'}(R \wedge \ddc \Phi)|\\
& \lesssim & \|R \wedge \ddc \Phi\|_{\Cc^{-2}}^\lambda \lesssim
\|R\|_{\Cc^{-2}}^\lambda.
\end{eqnarray*}
The result follows.
\endproof

Moderate currents and moderate measures were introduced in 
\cite{DinhNguyenSibony, DinhSibony1}. With respect to test d.s.h. functions, moderate
measures have the same regularity as the Lebesgue measure does. 

\begin{definition}\rm
A positive measure $\nu$ on $X$ is {\it moderate} if there are
constants $\lambda>0$ and $A>0$ such that
$$\langle\nu, e^{\lambda |\phi|}\rangle\leq A$$
for every d.s.h. function $\phi$ on $X$ such that
$\|\phi\|_\DSH\leq 1$. A measure is {\it moderate} if it is a
difference of moderate positive measures. A positive closed
$(p,p)$-current $S$ is {\it moderate} if its trace measure
$S\wedge\omega^{k-p}$ is moderate and a current in $\Dc_p$ is {\it
  moderate} if it is a difference of moderate positive closed
$(p,p)$-currents. 
\end{definition}

\begin{proposition} \label{prop_holder_moderate}
Let $S$ be a positive closed $(p,p)$-current. Assume that $S$ has a
H{\"o}lder continuous super-potential. Then $S$ is
moderate. 
\end{proposition}
\proof
By Proposition \ref{prop_holder_wedge}, the trace measure of $S$ has a
H{\"o}lder continuous super-potential. Replacing $S$ by its trace measure
allows to assume that $S$ is a
positive measure, i.e. $p=k$. Let $\phi$ be a d.s.h. function with
$\|\phi\|_\DSH\leq 1$. Define
$\psi_M:=\min(|\phi|,M)-\min(|\phi|,M-1)$ for $M\geq 1$. 
Observe that $0\leq\psi_M\leq 1$ and that $\psi_M=0$ on $\{|\phi|\leq
M-1\}$, also $\psi_M$ is larger than or equal
to the characteristic function $\rho_M$ of $\{|\phi|\geq M\}$.
Moreover, the DSH-norm of $\psi_M$ is bounded by a constant independent of
$\phi$ and $M$, see e.g. \cite{DinhSibony6}.
We want to prove that
$\langle S,e^{\lambda|\phi|}\rangle\leq A$ for some
positive constants $\lambda$ and $A$. 
So, it is enough to check that $\langle S,\rho_M\rangle\lesssim
e^{-\lambda M}$ for some (other) positive constant $\lambda$. For this
purpose, we will show that $\langle S,\psi_M\rangle\lesssim
e^{-\lambda M}$. 

By Proposition \ref{prop_exp}, the volume of the support of $\psi_M$ is
$\lesssim e^{-\lambda M}$ since it is contained in $\{|\phi|\geq
M-1\}$. Therefore, the estimate $\langle S,\psi_M\rangle\lesssim
e^{-\lambda M}$ is clear when $S$ is a form with bounded
$\|\cdot\|_\infty$-norm because $0\leq\psi_M\leq 1$. 
Subtracting from $S$ a smooth form allows to assume that $[S]=0$ but
we loose here the positivity of $S$. 
Recall that the super-potential $\Uc_S$ of $S$ is H{\"o}lder
continuous and that $\psi_M$ has a bounded DSH-norm. We have for some
constant $\lambda>0$
$$\langle S,\psi_M\rangle =\Uc_S(\ddc\psi_M)\lesssim
\|\ddc\psi_M\|_{\Cc^{-2}}^\lambda.$$
On the other hand, if $\Phi$ is a test form with $\|\Phi\|_{\Cc^2}\leq
1$ then
$$\|\ddc\psi_M\|_{\Cc^{-2}}=\sup_\Phi|\langle \ddc\psi_M,\Phi\rangle|
=\sup_\Phi|\langle \psi_M,\ddc \Phi\rangle|\lesssim e^{-\lambda M},$$
where the last inequality follows from the above volume estimate 
of the support of $\psi_M$.
This completes the proof.
\endproof

\begin{proposition} \label{prop_holder_hausdorff}
Let $S$ be a positive closed $(p,p)$-current with H{\"o}lder continuous
super-potentials on $X$. Assume that the manifold $X$ is projective. 
Then the Hausdorff dimension of $S$ is strictly
larger than $2(k-p)$. More precisely, the trace measure
$S\wedge\omega^{k-p}$ has no mass on sets of finite Hausdorff measure
of dimension $2(k-p)+\epsilon$ for $\epsilon>0$ small enough.
\end{proposition}

We will need the following lemma where we use that $X$ is projective.

\begin{lemma}
Let $A>0$ be a constant large enough and $r_0>0$ a constant small
enough. If  $B_{r_0}$, $B_r$ are concentric balls of radius $r_0$,
$r$ respectively, $r\ll r_0$, then there is 
a positive smooth form $\Phi$ of bidegree $(k-p,k-p)$ supported in $B_{r_0}$ with $\Phi\geq
\omega^{k-p}$ on $B_r$ and such that 
$$\|\Phi\| \leq A r^{2k-2p+2},\quad
\|\ddc\Phi\|_*\leq A r^{2k-2p}\quad \mbox{and}\quad \|\ddc\Phi\|_{\Cc^{-1}}\leq A r^{2k-2p+1}.$$ 
\end{lemma} 
\proof
The case where $X$ is the projective space $\P^k$ is proved in
\cite[Lemma 3.3.7]{DinhSibony10}. We will deduce the lemma from this
particular case. Since $r_0$ is small, we can choose a finite
family of holomorphic maps from $X$ onto $\P^k$ such that every ball
of radius $3r_0$ is sent injectively to $\P^k$ by at least one of these
maps. Let $\Pi:X\rightarrow \P^k$ be such a map corresponding to the
ball $B_{3r_0}$ with the same center as the
considered balls $B_r$ and $B_{r_0}$. Then, $\Pi(B_{r_0})$
contains a ball $B'$ in $\P^k$ of radius $\gtrsim r_0$ and
$\Pi(B_r)$ is contained in a ball $B''$ 
of radius $\lesssim r$. Let $\Psi$ be a form satisfying the lemma for
$\P^k$, $B'$, $B''$ and for a fixed K{\"a}hler metric on $\P^k$. The choice of
$\Pi$ implies that the jacobian of $(\Pi_{|B_{2r_0}})^{-1}$ is bounded from
below and from above by positive constants. Therefore, the form
$\Phi:=\Pi_{|B_{r_0}}^*(\Psi)$ is positive with support in
$B_{r_0}$. It satisfies $\Phi\gtrsim \omega^{k-p}$ on $B_r$ and $\|\Phi\| \lesssim r^{2k-2p+2}$,
$\|\ddc\Phi\|_*\lesssim r^{2k-2p}$ on $X$. Multiplying
$\Psi$ by a constant allows to get $\Phi\geq \omega^{k-p}$ on
$B_r$. Finally, it remains to check the inequality $\|\ddc\Phi\|_{\Cc^{-1}}\lesssim
r^{2k-2p+1}$. We have to show that
$\sup_\Omega|\langle \ddc\Phi,\Omega\rangle|\lesssim r^{2k-2p+1}$
for smooth test form $\Omega$ with
$\|\Omega\|_{\Cc^1}\leq 1$. Since $\Phi$
is supported in $B_{r_0}$, it is enough to consider $\Omega$ with
support in $B_{2r_0}$. In that case, the desired estimate is
deduced from the analogous estimate for $\Psi$ in $\P^k$.
\endproof

\noindent
{\bf End of the proof of Proposition \ref{prop_holder_hausdorff}.} 
Fix a constant $\epsilon>0$ small enough. We only have to prove that $\int_{B_r}
S\wedge\omega^{k-p}\lesssim r^{2k-2p+\epsilon}$ for $r$
small, see e.g. \cite{Sibony}. 
Using the previous lemma, it suffices to check that  
$\langle S,\Phi\rangle\lesssim r^{2k-2p+\epsilon}$.
The estimate is clear when $S$ is
smooth. Subtracting from $S$ a smooth form allows to assume that
$[S]=0$ but we loose the positivity of $S$.
Let $\Uc_S$ be the super-potential of $S$. Since $\Uc_S$ is H{\"o}lder
continuous and $\|r^{2p-2k}\ddc\Phi\|_*\leq A$, we have
\begin{eqnarray*}
\langle S,\Phi\rangle & = & \Uc_S(\ddc\Phi) =
r^{2k-2p}\Uc_S(r^{2p-2k}\ddc\Phi) \\
& \lesssim & r^{2k-2p}
\|r^{2p-2k}\ddc\Phi\|_{\Cc^{-1}}^\epsilon\lesssim r^{2k-2p+\epsilon}.
\end{eqnarray*}
This implies the proposition.
\hfill $\square$


\section{Dynamics of automorphisms} \label{section_dynamics}

In this section, we study the dynamics of automorphisms on compact K{\"a}hler manifolds.
The main dynamical objects (Green currents and equilibrium measure)
were constructed by the authors in \cite{DinhSibony4}. 
In \cite{Guedj}, under some extra hypothesis, Guedj gives another construction of the Green current
of some bidegree and of the equilibrium measure.
Here, the theory of
super-potentials allows us to obtain a new construction and to prove some fine
properties of these dynamical objects.

\subsection{Action on currents and on cohomology groups} \label{section_linear}

We first give some basic properties of linear maps. Their proofs are
left to the reader. Recall that a {\it Jordan block}
$J_{\lambda,m}$ is a square complex matrix
$(a_{ij})_{1\leq i,j\leq m}$ such that
$a_{ij}=\lambda$ if $i=j$, $a_{ij}=1$ if $j=i+1$ and $a_{ij}=0$ 
otherwise. 
If $\lambda\not=0$, the entry of index
$(1,m)$ of $J_{\lambda,m}^n$ is 
equal to ${n\choose m-1}\lambda^{n-m+1}$ when $n\geq m-1$. This is
the only entry of order
$n^{m-1}|\lambda|^n$,
the other ones have order at most 
equal to $n^{m-2}|\lambda|^n$.
We have 
$$\|J_{\lambda,m}^n\|\sim
{n \choose m-1}|\lambda|^{n-m+1} \sim n^{m-1}|\lambda|^n.$$ 
The eigenspace of $J_{\lambda,m}$
associated to the unique eigenvalue $\lambda$ is a complex line.

Consider a linear automorphism $L$ of a real space $E\simeq \R^h$. We assume there is an
open convex cone $\Kc$ in $E$ which is salient, i.e. $\overline\Kc\cap
-\overline\Kc=\{0\}$, and totally invariant by $L$, i.e. $L(\Kc)=\Kc$. For a
fixed basis of $E$, $L$ is associated to an invertible square matrix
with real coefficients. 
One can extend $L$ to an automorphism of
$E^\C:=E\otimes_\R\C\simeq\C^h$. Then, there is a complex basis of $E^\C$ such
that the associated matrix of $L$ is a Jordan matrix, i.e. a block
diagonal matrix whose blocks are Jordan blocks. In other words,
one can decompose $E^\C$ into direct sum of complex subspaces $E_l^\C$ which are
invariant by $L$:
$$E^\C=\bigoplus_{1\leq l\leq r} E^\C_l\quad \mbox{with}\quad \dim E^\C_l=m_l\quad
\mbox{and}\quad \sum_{l=1}^r m_l=h,$$
such that the restriction $L_l$ of $L$ to $E^\C_l$ is defined by
a Jordan block $J_{\lambda_l,m_l}$.

Up to a permutation of the $E^\C_l$, we can assume that 
the $(|\lambda_l|,m_l)$ are ordered so that
either $|\lambda_l|>|\lambda_{l+1}|$ or  $|\lambda_l|=|\lambda_{l+1}|$
and $m_l\geq m_{l+1}$ for every $1\leq l\leq r-1$.
 We say that $J_{\lambda_l,m_l}$
is a  {\it dominant} Jordan block if $(|\lambda_l|,m_l)= (|\lambda_1|,m_1)$ and in
that case we say that $\lambda_l$ 
is a {\it dominant eigenvalue} of $L$. 
Let $\nu$ be the integer such that
$J_{\lambda_1,m_1}$, $\ldots$, $J_{\lambda_\nu,m_\nu}$ are the dominant
Jordan blocks. 
The positive number 
$\lambda:=|\lambda_1|$ is the {\it spectral radius} of
$L$. The integer $m:=m_1$ is called the {\it multiplicity of
the spectral radius}. Since $L$ preserves the salient open cone $\Kc$,
$\lambda$ is a dominant eigenvalue of $L$. Moreover, the Perron-Frobenius
theorem implies that $L$ admits an eigenvector in $\overline\Kc$
associated to the dominant eigenvalue $\lambda$. It is clear that $\|L^n\|\sim
n^{m_1-1}\lambda^n$.
Let $\widetilde E_l^\C$ be the hyperplane generated by the
first $m_l-1$ vectors of the basis of $E_l^\C$ associated to the Jordan
form. We have $\|L^n v\|\sim
n^{m_1-1}\lambda^n$
for any vector $v\not\in
\widetilde E_1^\C\oplus \cdots\oplus \widetilde E_\nu^\C\oplus 
E^\C_{\nu+1}\oplus\cdots \oplus E^\C_r$, in particular for $v\in\Kc$.

Let $F_l^\C$ denote the eigenspace of
$L_l$ which is a complex line. We say that $F^\C:=F^\C_1\oplus \cdots
\oplus F^\C_\nu$ is the {\it dominant eigenspace}. Define
$H^\C:=\oplus F^\C_l$ with $1\leq l\leq \nu$ and $\lambda_l=\lambda$. This is the {\it 
strictly dominant eigenspace} of $L$. Define also $F:=F^\C\cap E$ and $H:=H^\C\cap E$. One can check that
$F^\C=F\otimes_\R\C$ and $H^\C=H\otimes_R\C$. 
The previous spaces are 
invariant under $L$. 

For any $1\leq l\leq \nu$, there is a unique 
$\theta_l\in \mathbb{S}:=\R/2\pi\Z$ such that
$\lambda_l=\lambda\exp(i\theta_l)$. We say that 
$\theta:=(\theta_1,\ldots,\theta_\nu)\in \mathbb{S}^\nu$
is the {\it dominant direction} of $L$.
The dominant direction of $L^n$ is equal to $n\theta$. 
Denote by $\Theta$ the closed subgroup
of $\mathbb{S}^\nu$ generated by $\theta$.
It is a finite union
of real tori.
The orbit of each point
$\theta'\in \Theta$ under the translation
$\theta'\mapsto \theta'+\theta$ is dense in $\Theta$.  
If $\lambda_l=\lambda$ for every $1\leq l\leq\nu$,
we have $F^\C=H^\C$, $\theta=0$ and $\Theta=\{0\}$.
Define 
$$\widehat L_N:=\frac{1}{N}\sum_{n=1}^N {L^n\over n^{m-1}\lambda^n}.$$
We have the following proposition, see also \cite{DinhSibony4}.

\begin{proposition} \label{prop_iterate_linear}
The sequence $(\widehat L_N)$ converges to 
a surjective real linear map $\widehat
L_\infty:E\rightarrow H$.
Let $(n_i)$ be an increasing sequence of integers. Then
$(n_i^{1-m}\lambda^{-n_i}L^{n_i})$ converges if and only if
$(n_i\theta)$ converges. Moreover, any limit of
$(n^{1-m}\lambda^{-n}L^n)$ is a 
surjective real linear map from $E$ to $F$. 
\end{proposition}

Note that surjective linear maps are open and the image of $\Kc$ by
such a map is an open convex cone.

We will apply the previous result to the action of a holomorphic map on
cohomology groups.
Let $f$ be an automorphism on a compact K{\"a}hler manifold $(X,\omega)$
of dimension $k$. The pull-back operator $f^*$ acts on smooth forms
and on currents. It
commutes with $\partial$, $\dbar$ and preserves positivity. Therefore, $f^*$ acts as a
linear automorphism on
 $H^{q,q}(X,\R)$. The operator push-forward $f_*$
is defined in the same way. It coincides with the pull-back
$(f^{-1})^*$ by
$f^{-1}$. 

Recall that {\it the dynamical degree of order $q$} of $f$ is the
spectral radius of $f^*$ acting on $H^{q,q}(X,\R)$. Let us denote by
$d_q(f)$ (or $d_q$ if there is no confusion) this degree. We have
$d_q(f^n)=d_q(f)^n$ for $n\geq 1$ and $d_0(f)=d_k(f)=1$. 
An inequality due to Khovanskii, Teissier and Gromov \cite{Gromov1} implies
that the function $q\mapsto \log d_q$ is concave on $0\leq q\leq k$,
see also \cite{Guedj}. 
In particular, there are two integers $p$ and $p'$ with $0\leq p\leq p'\leq k$ such that
$$1=d_0<\cdots< d_p=\cdots=d_{p'}>\cdots>d_k=1.$$
By Gromov and Yomdin
\cite{Gromov2, Yomdin}, the dynamical degrees are related to the topological entropy
$h_t(f)$ of $f$ by the formula $h_t(f)=\max_q \log d_q$, see also
\cite{DinhSibony3} for a more general context.

Let $\Kc$ be the convex cone of the classes in $H^{q,q}(X,\R)$
associated to strictly positive closed $(q,q)$-forms. Then $\Kc$ is salient and
totally invariant by $f^*$. Hence, we can apply Proposition \ref{prop_iterate_linear} to
$f^*$ (one can also apply it to the
cone of the classes associated to strictly positive closed $(q,q)$-currents).
Recall that the bilinear form $\smile$ on $H^{q,q}(X,\R)\times
H^{k-q,k-q}(X,\R)$ given by
$$([\beta],[\beta'])\mapsto [\beta]\smile[\beta']:=\int_X\beta\wedge \beta'$$
is non-degenerate. Moreover, we have
$f^*[\beta]\smile[\beta']=[\beta]\smile f_*[\beta']$. So, if we consider two
basis of  $H^{q,q}(X,\R)$ and of
$H^{k-q,k-q}(X,\R)$ which are dual with respect to $\smile$, then $f^*$
acting on $H^{q,q}(X,\R)$ and $f_*$
acting on $H^{k-q,k-q}(X,\R)$ are given by the same matrix. Therefore,
these operators have the same spectral radius $d_q(f)=d_{k-q}(f^{-1})$
with the same multiplicity $m$. 

\begin{lemma} \label{lemma_mass_current}
If $S$ is a current in $\Dc_q$ then 
$$\|(f^n)^*S\|_*\leq \kappa n^{m-1}d_q^n\|S\|_*$$
where $\kappa>0$ is a constant independent of $S$. Moreover, if $S$ is
a strictly positive current, we have
$\|(f^n)^*(S)\|\sim n^{m-1}d_q^n$.
\end{lemma}
\proof
We can assume that $S$ is positive. The mass of a positive
closed current can be computed cohomologically. Therefore,
$$\|(f^n)^*S\| = (f^n)^*[S]\smile [\omega^{k-q}] \lesssim
\|(f^n)^*\|\|[S]\|\lesssim n^{m-1}d_q^n\|S\|.$$
This gives the first part of the lemma. For the second one, if $S$ is
strictly positive then $[S]$ is in the interior of the cone of the
classes of positive closed currents. We deduce from the above discussion on the
linear operator $L$ that $\|(f^n)^*[S]\|\sim n^{m-1}d_q^n$. The result
follows.
\endproof

Note that the previous lemma allows to compute the dynamical degrees
using the following formula
$$d_q(f)=\lim_{n\rightarrow\infty} \Big[\int_X (f^n)^*\omega^q\wedge
\omega^{k-q}\Big]^{1/n} =\lim_{n\rightarrow\infty} \Big[\int_X \omega^q\wedge
(f^n)_*\omega^{k-q}\Big]^{1/n}.$$

\subsection{Construction of Green currents} \label{section_green}

In this section, we give a new construction of the Green currents using the
super-potentials. This approach permits to establish some new properties of the
Green currents. The result can be extended to open non-invertible
maps, see \cite{DinhSibony8} for the pull-back operator on currents by
non-invertible maps. We use here the notations introduced in Section \ref{section_linear}.

\begin{theorem} \label{th_green_construction}
Let $f$ be a holomorphic automorphism on a compact K{\"a}hler manifold $(X,\omega)$. 
Let $d_s$ be the dynamical degrees of $f$ and $q$ an
integer such that $d_{q-1}<d_q$. Let $F$ (resp. $H$) denote the real 
dominant (resp. strictly dominant) subspace associated
to the operator $f^*$ on $H^{q,q}(X,\R)$. Then each class $c$ of $F$
can be represented by a current $T_c$ in $\Dc_q$ with H{\"o}lder continuous
super-potentials which depends linearly on
$c$ and satisfies $f^*(T_c)=T_{f^*(c)}$. In particular,
we have $f^*(T_c)=d_qT_c$ for $c\in H$. The set of the classes $c$ in
$F$ (resp. in $H$) with $T_c$ positive is a closed convex cone with non-empty interior.
\end{theorem}

The cone of positive closed
currents $T_c$ with $c\in H$ is a closed cone of finite
dimension. We say that $T_c$  is {\it a Green current of order
  $q$} of $f$ if $T_c$ is a non-zero positif current (this implies
that $c\not=0$).
By Proposition \ref{prop_holder_moderate}, the Green currents are moderate.
We will see that Green currents are the only positive closed
currents in their cohomology classes.

Consider now the action of $f^*$ on $H^{q,q}(X,\R)$ as described in
Section \ref{section_linear}. Let  $m$ denote the
multiplicity of its spectral radius $d_q$.

\begin{proposition} \label{prop_green_construction}
Let $S$ be a current in $\Dc_q$ with a continuous super-potential. Let
$(n_i)$ be an increasing sequence of integers. 
Assume that $n_i^{1-m}d_q^{-n_i}(f^{n_i})^*[S]$ converge to some class
$c$ in $H^{q,q}(X,\R)$.  
Then $n_i^{1-m}d_q^{-n_i}(f^{n_i})^*(S)$ converge SP-uniformly to a current $T_c$ in $\Dc_q$ which
depends only on $c$.
\end{proposition}

Let $\alpha=\{\alpha_1,\ldots,\alpha_h\}$ 
be a family of smooth closed real $(q,q)$-forms such that
$[\alpha]=\{[\alpha_1],\ldots,[\alpha_h]\}$ is a basis of
$H^{q,q}(X,\R)$ where $h$ is the dimension of $H^{q,q}(X,\R)$. 
In what follows, we consider the super-potentials normalized
by $\alpha$ as in Section \ref{section_sp}. Let $M$ denote the
$h\times h$ matrix whose column of index $j$ is given by the
coordinates of $f^*[\alpha_j]$ with respect to the basis $[\alpha]$. 
Let $\Uc_j$ denote the super-potential of $f^*(\alpha_j)$ and define
$\Uc:=(\Uc_1,\ldots,\Uc_h)$. Let
$A=\trans(a_1,\ldots,a_h)$ denote the coordinates of $[S]$ in the basis
$[\alpha]$ and $\Uc_S$, $\Uc_{S_n}$ the super-potentials of $S$ and
of $S_n:=(f^n)^*S$ respectively. Denote also by $\Lambda$ the operator $f_*$ acting on $\Dc^0_{k-q+1}$.

\begin{lemma} \label{lemma_sp_iterate}
We have
$$\Uc_{S_n}=\sum_{l=0}^{n-1} (\Uc\circ \Lambda^l) M^{n-l-1} A  + \Uc_S\circ \Lambda^n.$$
\end{lemma}
\proof
The proof is by induction. For $n=0$, we have $S_0=S$ and the lemma is
clear. Assume the lemma for $n$. We show it for $n+1$. Let $R$ be a
smooth form in $\Dc^0_{k-q+1}$ and $U$ a smooth potential of $R$ normalized
by $\alpha$. So, $f_*(U)$ is a potential of $\Lambda(R)$ but it is not
normalized. Let $\alpha^\vee=\{\alpha_1^\vee,\ldots,\alpha_h^\vee\}$
be a family of smooth real closed forms such that $[\alpha^\vee]$ is the basis of $H^{k-q,k-q}(X,\R)$
which is dual to $[\alpha]$ with respect to $\smile$. 
Define 
$$b_j:=\langle \alpha_j,f_*(U)\rangle = \langle f^*(\alpha_j),
U\rangle =\Uc_j(R)$$
and 
$$b:=\trans (b_1,\ldots,b_m)= \trans \Uc(R).$$ 
Then
$U':=f_*(U)-\alpha^\vee b$ is a potential of $\Lambda(R)$ normalized by
$\alpha$. We obtain using the induction hypothesis 
\begin{eqnarray*}
\Uc_{S_{n+1}}(R) & = & \langle S_{n+1},U\rangle = \langle
f^*(S_n),U\rangle = \langle S_n, f_*(U)\rangle \\
& = & \langle S_n,U'\rangle +\langle S_n, \alpha^\vee b\rangle
=\Uc_{S_n}(\Lambda(R)) +\langle S_n, \alpha^\vee b\rangle\\
& = & \sum_{l=1}^n \Uc(\Lambda^l(R))M^{n-l} A + \Uc_S(\Lambda^{n+1}(R)) +\langle S_n, \alpha^\vee b\rangle.
\end{eqnarray*}
We only have to check that the last integral satisfies
$$\langle S_n, \alpha^\vee b\rangle= \Uc(R)M^nA.$$
Observe that the  integral $\langle S_n, \alpha^\vee b\rangle$ can be
computed cohomologically. Since $S$ is cohomologous to $\alpha A$, by
definition of $M$, $S_n=(f^n)^*S$ is cohomologous to
$\alpha M^nA$. Hence 
$$\langle S_n, \alpha^\vee b\rangle=\langle \alpha M^nA,
\alpha^\vee b\rangle=\trans b M^nA=\Uc(R)M^nA.$$
This completes the proof.
\endproof

\noindent
{\bf End of the proof of Proposition \ref{prop_green_construction}.} 
By Proposition \ref{prop_iterate_linear} the limit $c$ of $n_i^{1-m}d_q^{-n_i} (f^{n_i})[S]$
is a class in $F$. Write $c=\trans (c_1,\ldots, c_h)$ with respect to
the basis $[\alpha]$. Then, $n_i^{1-m}d_q^{-n_i} M^{n_i} A$ converges to
$c$. 

By Lemma \ref{lemma_sp_iterate}, the super-potential $\Uc_{n_i}$ of
$n_i^{1-m}d_q^{-n_i} (f^{n_i})^*(S)$ is equal to
\begin{eqnarray}
\Uc_{n_i} & = & n_i^{1-m}d_q^{-n_i}
\Big[\sum_{l=0}^{n_i-1} (\Uc\circ \Lambda^l)  M^{n_i-l-1} A + \Uc_S\circ
\Lambda^{n_i}\Big] \nonumber\\
& = & 
\sum_{l=0}^{n_i-1} (\Uc\circ \Lambda^l) {M^{n_i-l-1}A\over n_i^{m-1}d_q^{n_i}}
+ n_i^{1-m}d_q^{-n_i} \Uc_S\circ \Lambda^{n_i} \label{eq_sp}.
\end{eqnarray}
Since $\Uc_S$ is continuous, we have 
\begin{equation}
|\Uc_S\circ \Lambda^n(R)|\lesssim\|\Lambda^n(R)\|_*\lesssim
\delta^n\|R\|_*, \label{eq_sp_iterate}
\end{equation}
where $d_{q-1}<\delta<d_q$ is a fixed constant.
It follows that the last term in (\ref{eq_sp}) tends uniformly to 0 on
$\ast$-bounded sets of $R$.

Recall that $\|M^n\|\sim n^{m-1}d_q^n$. Analogous estimates as in (\ref{eq_sp_iterate}) for
 $\Uc_j$ imply that
\begin{equation}
\Big\| {M^{n_i-l-1}A\over n_i^{m-1}d_q^{n_i}}\Big\|
\lesssim d_q^{-l} \quad \mbox{and}\quad 
\Big|(\Uc_j\circ \Lambda^l) {M^{n_i-l-1}A\over n_i^{m-1}d_q^{n_i}}\Big|
\lesssim \delta^l d_q^{-l}. \label{eq_matrix_iterate}
\end{equation}
Since $\sum_{l\geq 0} \delta^ld_q^{-l}$ converges, 
we can apply the Lebesgue convergence theorem for the sum in
(\ref{eq_sp}). We obtain the uniform convergence on $\ast$-bounded sets: 
$$\lim_{i\rightarrow\infty} \Uc_{n_i}=\sum_{l\geq 0}  (\Uc\circ
\Lambda^l) M^{-l-1}c.$$
The last series converge because (\ref{eq_matrix_iterate}) implies that
$\|M^{-l-1}c\|\lesssim d_q^{-l}$. One can also obtain this inequality using the fact
that $c$ is a vector in $F$ and that the matrix of $f^*_{|F}$ 
is conjugate to a diagonal matrix whose eigenvalues are of
modulus $d_q$. 

Hence, the sequence $(n_i^{1-m} d_q^{-n_i} (f^{n_i})^*(S))$ converges
to some current $T_c$. Moreover, the last series defines a
super-potential $\Uc_{T_c}$ of $T_c$ and $\Uc_{n_i}$ converge
  to $\Uc_{T_c}$ uniformly on $\ast$-bounded sets of $R$. 
Hence, the  convergence of  $(n_i^{1-m} d_q^{-n_i} (f^{n_i})^*(S))$  is
  SP-uniform. Since $\Uc_{T_c}$ depends only on the class $c$, by
  Proposition \ref{prop_current_sp},
$T_c$ depends only on this class.
\hfill $\square$

\begin{proposition} \label{prop_green_inverse}
Let $(n_i)$ be an increasing sequence such that
  $(n_i^{1-m}d_q^{-n_i}(f^{n_i})^*)$ converges on
  $H^{q,q}(X,\R)$. Then for any class $c\in F$ there is a smooth form
  $S$ in $\Dc_q$ such that $n_i^{1-m}d_q^{-n_i}(f^{n_i})^*(S)$
  converge SP-uniformly to $T_c$.
\end{proposition}
\proof
By Proposition \ref{prop_iterate_linear}, 
one can find a smooth form $S$ in $\Dc_q$ such that
$n_i^{1-m}d_q^{-n_i}(f^{n_i})^*[S]$ converge to $c$. It is enough to
apply Proposition \ref{prop_green_construction}.
\endproof

\begin{lemma} \label{lemma_holder}
The current $T_c$ in Proposition \ref{prop_green_construction} has H{\"o}lder continuous potential.
\end{lemma}
\proof
We follow the approach in \cite{DinhSibony4} and \cite{DinhSibony10}.
Let $R$ be a smooth form in $\Dc^0_{k-q+1}$ such that $\|R\|_*\leq 1$.
Since $f^*(\alpha_j)$ is smooth, its super-potential $\Uc_j$ is
Lipschitz in the sense that $\|\Uc_j(R)\|\leq \kappa \|R\|_{\Cc^{-1}}$ with
$\kappa>1$ independent of $R$. By definition of $\|\cdot \|_{\Cc^{-1}}$, we also
  have $\|\Lambda(R)\|_{\Cc^{-1}}\leq \kappa\|R\|_{\Cc^{-1}}$ for some
    constant $\kappa>1$. 

Let $\delta$ be a constant as above with $d_{q-1}<\delta<d_q$. Define
$\rho:=\delta d_q^{-1}$, 
$\widetilde\Lambda:=\delta^{-1}\Lambda$, $\lambda:=-\log\rho(\log \kappa-\log\rho)^{-1}$ and $N_0$ the integer
part of $(\lambda-1)\log\|R\|_{\Cc^{-1}}(\log \kappa)^{-1}$. Then the sequence
$(\widetilde\Lambda^l)_{l\geq 0}$ is bounded with respect to the $\|\cdot\|_*$-norm. 
For $\|R\|_{\Cc^{-1}}$ small enough, we have since
$\|M^{-l-1}c\|\lesssim d_q^{-l}$
\begin{eqnarray*}|\sum_{l\geq 0}
(\Uc\circ \Lambda^l(R)) M^{-l-1}c | & \lesssim &  \sum_{l=0}^{N_0}\rho^l
\|\Uc\circ \widetilde \Lambda^l(R)\|+\sum_{l>N_0}
\rho^l\|\Uc\circ \widetilde\Lambda^l(R)\|\\
&\lesssim&
\Big(\sum_{l=0}^{N_0}\rho^l\kappa^l\Big)\|R\|_{\Cc^{-1}}+\sum_{l>N_0} \rho^l\\
&\lesssim & \kappa^{N_0} \|R\|_{\Cc^{-1}}+\rho^{N_0}\lesssim \|R\|_{\Cc^{-1}}^\lambda.
\end{eqnarray*}
Therefore, $\sum_{l\geq 0} (\Uc\circ \Lambda^l)M^{-l-1}c$ is
$\lambda$-H{\"o}lder continuous with respect to $\dist_{-1}$. 
\endproof

\noindent
{\bf End of the proof of Theorem \ref{th_green_construction}.} 
By Proposition \ref{prop_green_inverse}, $T_c$
depends linearly on $c$ because it depends linearly on $S$. 
By Lemma \ref{lemma_holder}, $T_c$ has
H{\"o}lder continuous super-potentials. 
Observe that if $n_i^{1-m}d_q^{-n_i} (f^{n_i})^*[S]$ converge to $c$
then $n_i^{1-m}d_q^{-n_i} (f^{n_i})^*[f^*(S)]$ converge to $f^*(c)$.
Applying Proposition \ref{prop_green_construction} to
$S$ and to $f^*(S)$ yields
$$f^*(T_c)  =  f^*\Big(\lim_{i\rightarrow\infty} n_i^{1-m}d_q^{-n_i} (f^{n_i})^*(S)\Big) = 
\lim_{i\rightarrow\infty} n_i^{1-m}d_q^{-n_i} (f^{n_i})^*\big(f^*(S)\big) =T_{f^*(c)}.$$
If $c$ is in $H$, we have $f^*(c)=d_qc$. Hence, $f^*(T_c)=d_qT_c$. 

We deduce easily from the linear dependence of $T_c$ on $c$ that the
cone $\Cc_F$ (resp. $\Cc_H$) of the classes $c$ in $F$ (resp. in $H$) with $T_c$ positive is
convex and closed. It remains to prove that they have non-empty
interior. Observe that $\Cc_F$ contains the classes $c$ associated to $S$ smooth 
strictly positive and that the cone $\Kc$ in $H^{q,q}(X,\R)$ of
such forms $S$ is open. By Proposition \ref{prop_iterate_linear}, any limit $L_\infty$ of
$(n^{1-m}d_q^{-n}(f^n)^*)$ is an open map from $H^{q,q}(X,\R)$ to $F$.
Hence, $\Cc_F$ contains the cone $L_\infty(\Kc)$ which is open
in $F$. 

Consider as in Section \ref{section_linear} the operators
$$\widehat L_N:={1\over N}\sum_{n=1}^N n^{1-m}d_q^{-n} (f^n)^*$$
on $H^{q,q}(X,\R)$
and $\widehat L_\infty$ the limit of this sequence which is an open
map from $H^{q,q}(X,\R)$ to $H$. Observe that any class in $\widehat
L_\infty(\Kc)$ belongs to the closed convex cone generated by the
$L_\infty(\Kc)$. Hence, $\Cc_F$ contains $\widehat
L_\infty(\Kc)$. We deduce that $\Cc_H$, which is equal to
$\Cc_F\cap H$, contains the open cone $\widehat L_\infty(\Kc)$ of
$H$. This completes the proof of  Theorem \ref{th_green_construction}.
\hfill $\square$

\begin{proposition} \label{prop_cesaro_cv}
Let $S$ be a current in $\Dc_q$ with continuous super-potentials. Let
$\widehat L_N$ be as above and $\widehat c\in H$ the limit of the
sequence $(\widehat L_N[S])$. Then $\widehat L_N(S)$ converge
SP-uniformly to the current $T_{\widehat c}$. Moreover, any current
$T_{\widehat c}$ with $\widehat c\in H$ can be obtained as the limit
of  $(\widehat L_N[S])$ for some $S$ smooth in $\Dc_q$.
\end{proposition}
\proof
Proposition \ref{prop_green_construction} implies that any limit $T$ of $\widehat L_N(S)$
belongs to the space generated by the $T_c$ with $c\in F$. Hence, $T$
is equal to one of the current $T_c$. On the other hand, $T$ is a
current in the class $\widehat c$. We deduce that $c=\widehat c$ and
that $\widehat L_N(S)$ converge to $T_{\widehat c}$. The main point
here is to show that the convergence is SP-uniform. We follows 
the proof of Proposition  \ref{prop_green_construction}.

By Lemma \ref{lemma_sp_iterate} the super-potential $\Uc_{\widehat L_N(S)}$ of
$\widehat L_N(S)$ is equal to
\begin{eqnarray}
\Uc_{\widehat L_N(S)} & = & {1\over N} \sum_{n=1}^N n^{1-m}d_q^{-n}
\Big[\sum_{l=0}^{n-1} (\Uc\circ \Lambda^l)  M^{n-l-1} A + \Uc_S\circ
\Lambda^n\Big] \nonumber\\
& = & \sum_{l=0}^{N-1} d_q^{-l-1}(\Uc\circ \Lambda^l)\Big[{1\over N}\sum_{n=0}^{N-l-1}
(n+l+1)^{1-m}d_q^{-n}M^nA\Big] + \label{eq_cesaro_1}\\
&  & 
+  {1\over N} \sum_{n=1}^N n^{1-m}d_q^{-n}\Uc_S\circ
\Lambda^n. \label{eq_cesaro_2}
\end{eqnarray}
Since $\Uc_S$ is continuous, the quantity in (\ref{eq_cesaro_2}) tends
to 0 uniformly on $\ast$-bounded sets. By Proposition
\ref{prop_iterate_linear}, the term in the 
brackets in
(\ref{eq_cesaro_1}) converges to the vector of coordinates equal to
the coordinates of
$\widehat c$ in the basis $[\alpha]$. Denote also by $\widehat
c$ this vector. We deduce that the expression
in  (\ref{eq_cesaro_1}) converges uniformly on $\ast$-bounded sets to 
$$\Uc_{T_{\widehat c}}:=\sum_{l\geq 0} d_q^{-l-1} (\Uc\circ
\Lambda^l)\widehat c,$$
which defines a super-potential of $T_{\widehat c}$. Hence, the
convergence of $\widehat L_N(S)$ is SP-uniform.

The last assertion of the proposition is deduced from the surjectivity
of the map $\widehat L_\infty$ in Proposition \ref{prop_iterate_linear}.
\endproof


\subsection{Uniqueness of Green currents and equidistribution}

In this section, we will prove the  uniqueness of the Green currents
in their cohomology classes. 
We have the following general result.

\begin{theorem} \label{th_unique_green}
Let $f$ be a holomorphic automorphism of a compact
  K{\"a}hler manifold $(X,\omega)$ and $d_s$ 
  the dynamical degrees of $f$. Let $V$ be a subspace of
  $H^{q,q}(X,\R)$  invariant by $f^*$. Assume that $d_q>d_{q-1}$ and
  that all the (real and complex) eigenvalues
  of the restriction of $f^*$ to $V$ are
  of  modulus strictly larger than $d_{q-1}$.
Then each class in $V$ contains at  most one positive closed $(q,q)$-current.
\end{theorem}
\proof
Let $S$ and $S'$ be positive closed currents in the same class in $V$. Define 
$\lambda_n:=\|(f^n)_*(S)\|^{-1}=\|(f^n)_*(S')\|^{-1}$,
$S_n:= \lambda_n (f^n)_*(S)$ and $S_n':= \lambda_n (f^n)_*(S')$. 
The currents $S_n$, $S_n'$ are of mass 1. We have
$S=\lambda_n^{-1}(f^n)^*(S_n)$ and 
$S'=\lambda_n^{-1}(f^n)^*(S'_n)$. 
Let $\delta_1>\delta_2 >d_{q-1}$ be constants such that the 
eigenvalues of $f^*_{|V}$ have modulus strictly larger than
$\delta_1$. 
We deduce that the eigenvalues of $f_{*|V}$ have modulus strictly
smaller than $\delta_1^{-1}$. Therefore, Lemma
\ref{lemma_mass_current} applied to $f_*$ implies that 
$\lambda_n\gtrsim
\delta_1^n$. Let $\Uc_S$, $\Uc_{S'}$,
$\Uc_{S_n}$, $\Uc_{S_n'}$ denote the super-potentials of $S$, $S'$,
$S_n$ and $S_n'$ respectively.

Assume that $S\not=S'$. Proposition \ref{prop_current_sp} implies that
$\Uc_S\not=\Uc_{S'}$. Then, there is a smooth
form $R$ in $\Dc^0_{k-q+1}$ such that $\Uc_S(R)-\Uc_{S'}(R)\not=0$. If
we multiply $R$ by a constant, we can assume that   $\Uc_S(R)-\Uc_{S'}(R)=1$.
Since $S$ and $S'$ are cohomologous, they have the same coordinates
$A$ in the basis $[\alpha]$.  By Lemma \ref{lemma_sp_iterate}, we have
\begin{eqnarray*}
\Uc_{S_n}((f^n)_*R)-\Uc_{S_n'}((f^n)_*R) & = &
\Uc_{(f^n)^*S_n}(R)-\Uc_{(f^n)^*S_n'}(R) \\
& = & \lambda_n \Uc_S(R)-\lambda_n \Uc_{S'}(R)\\
& = & \lambda_n.
\end{eqnarray*}
Define $R_n:=\gamma_n^{-1} (f^n)_*(R)$ where $\gamma_n$ is the norm of
$(f^n)^*$ acting on $H^{q-1,q-1}(X,\R)$. 
Recall that $\gamma_n$ is also the norm of $(f^n)_*$ acting on
$H^{k-q+1,k-q+1}(X,\R)$. 
We have $\lim_{n\rightarrow\infty} \gamma_n^{1/n}=d_{q-1} <\delta_2$.
Moreover, $\|R_n\|_*$ is bounded by a constant independent of $n$. 
We have for $n$ large enough 
$$\Uc_{S_n}(R_n)-\Uc_{S'_n}(R_n)=\lambda_n\gamma_n^{-1}\geq 2\delta_1^n\delta_2^{-n}.$$
It follows that either $|\Uc_{S_n}(R_n)|\geq \delta_1^n\delta_2^{-n}$
or  $|\Uc_{S'_n}(R_n)|\geq \delta_1^n\delta_2^{-n}$.
On the other hand, if $\kappa$ is a fixed constant large enough, we have
$\|f_*(R)\|_{\Cc^1}\leq \kappa \|R\|_{\Cc^1}$ since $f$ is an automorphism, and by induction
$$\|R_n\|_{\Cc^1} = \gamma_n^{-1} \|(f^n)_*(R)\|_{\Cc^1}\lesssim
\widetilde \kappa^n,$$
for some constant $\widetilde \kappa$.
Theorem \ref{th_main_estimate}, applied to $S_n$
and $S_n'$, implies that
$\delta_1^n\delta_2^{-n}\lesssim n$. This is a contradiction because
$\delta_1>\delta_2$. 
\endproof

Observe that in Theorem \ref{th_unique_green}, by linearity, each class of $V$ contains at most one
current $T=T^+-T^-$ with $T^+$, $T^-$ positive closed and $[T^+]$,
$[T^-]$ in $V$. We apply this theorem to $V$ the maximal subspace of
$H^{q,q}(X,\R)$ where the eigenvalues of $f^*$ are of modulus
$d_q$. We obtain the following corollaries.

\begin{corollary} \label{cor_unique_green}
Let $f$, $d_s$ and $q$ be as in Theorem \ref{th_green_construction}. Then the Green
$(q,q)$-currents of $f$ are the unique positive closed currents 
in their cohomology classes. They are the only non-zero positive closed
$(q,q)$-currents which are invariant by $d_q^{-1}f^*$, i.e. satisfying the
equation  $d_q^{-1}f^*(T)=T$. 
\end{corollary}

\begin{corollary} Let $f$, $d_s$ and $q$ be as in Theorem
  \ref{th_green_construction}. Let $T$ be a Green $(q,q)$-current of
  $f$ and $\Cc_T$ the set of positive closed $(q,q)$-currents $S$
  such that $S\leq c T$ for some constant $c>0$. Then $\Cc_T$ is a
  salient convex closed cone of finite dimension. 
Moreover, each current in $\Cc_T$ is the unique positive
  closed current in its cohomology class. 
\end{corollary}
\proof
It is clear that $\Cc_T$ is a convex cone. It is salient since the
cone of of all positive closed $(q,q)$-currents is salient. Let $E^+$
denote the cone of the classes of currents $S$ in $\Cc_T$ and $E$ the
space generated by $E^+$. Then $E^+$ is convex and salient since it is contained in the cone
of the classes of positive closed currents. Since $T$ is a Green
current, it is invariant by $d_q^{-1} f^*$ and $d_q f_*$. 
If $v$ is a vector in $E^+$, then by definition of $E^+$,
$\|(f^n)_*v\|\lesssim d_q^{-n}$. Therefore, the eigenvalues of $f_*$
restricted to $E$ are of modulus at most equal to $d_q^{-1}$. We
deduce that
all eigenvalues of $f^*$ have
modulus equal to $d_q$. By Theorem \ref{th_unique_green}, $S$ is the only positive
closed current in $[S]$. Moreover, in $\overline E^+$, the current $S$
depends linearly on its class. We deduce from the correspondence
$S\leftrightarrow [S]$ and the definition of $E^+$  that $\overline
E^+=E^+$. So, $\Cc_T$ is closed. 
\endproof

The following results can be applied to the currents of integration on
subvarieties of pure codimension $q$ of $X$, and give equidistribution
properties of their images by $f^{-n}$.

\begin{corollary} \label{corollary_equi_1}
Let $f$, $d_s$ and $q$  be as  in Theorem
\ref{th_green_construction}. Let $m$ denote the multiplicity of the
spectral radius of $f^*$ on $H^{q,q}(X,\R)$. Let $(S_i)$ be a sequence of
positive closed $(q,q)$-currents.
If $(n_i)$ is an increasing
sequence of integers such that $n_i^{1-m}d_q^{-n_i}(f^{n_i})^*[S_i]$
converge in $H^{q,q}(X,\R)$, 
then $n_i^{1-m}d_q^{-n_i}(f^{n_i})^*(S_i)$ converge either to $0$ or to a Green
$(q,q)$-current.
\end{corollary}
\proof
Let $c$ denote the limit of $n_i^{1-m} d_p^{-n_i}
(f^{n_i})^*[S_i]$. Observe that the sequence of $n_i^{1-m} d_p^{-n_i}
(f^{n_i})^*[S_i]$ is bounded. Hence,
the currents $n_i^{1-m}d_q^{-n_i}(f^{n_i})^*(S_i)$
have mass bounded by a constant independent of $i$. Then, we can
extract convergence subsequences. All the limit currents are in the
same class $c$ of $F$. By
Theorem \ref{th_unique_green}, they are equal. This implies the result.
\endproof

\begin{corollary} \label{cor_cesaro_unique}
Let $f$, $d_s$, $q$ and $m$ be as in Corollary 
\ref{corollary_equi_1}. 
Let $S$ be a positive closed $(q,q)$-current on $X$. Then the sequence
$${1\over N} \sum_{n=1}^N n^{1-m} d_q^{-n}(f^n)^*(S)$$
converges either to $0$ or to a Green $(q,q)$-current of $f$.
\end{corollary}
\proof
It is enough to apply Theorem \ref{th_unique_green} and to observe 
that by Proposition \ref{prop_iterate_linear}, the sequence 
$${1\over N} \sum_{n=1}^N n^{1-m} d_q^{-n}(f^n)^*[S]$$
converges to a class in $H$. Note that when all the dominant
eigenvalues of $f^*$ on $H^{q,q}(X,\R)$ are equal to $d_q$ (i.e. real
positive), then $n^{1-m}d_q^{-n}(f^n)^*[S]$ converges to a class $c$ in
$H$. Therefore,
$n^{1-m}d_q^{-n}(f^n)^*(S)$ converge to 0 if $c=0$ or to a Green current otherwise.
\endproof


\subsection{Equilibrium measure, mixing and hyperbolicity}

In this section, we assume that $f$ admits a dynamical degree $d_p$ strictly
larger than the other ones. We have
$$1=d_0<\cdots<d_p>\cdots>d_k=1.$$
Then, we can construct Green $(q,q)$-currents of $f$ for $1\leq q\leq
p$ and Green $(q,q)$-currents associated to
$f^{-1}$ for $1\leq q\leq k-p$. 

If $T^+$ is a Green
$(p,p)$-current of $f$ and $T^-$ is a Green $(k-p,k-p)$-current
associated to $f^{-1}$, then as it is noticed in \cite{DinhSibony4}, we can define the
intersection $T^+\wedge T^-$, see also Section
\ref{section_intersection}. 
This gives an invariant measure. However, we cannot prove that this measure does not vanish.
We introduced in \cite{DinhSibony4} another construction which always
gives an ergodic probability measure\footnote{there is a slip at the
  end of \cite[p.310]{DinhSibony4}; the measure that we constructed is only
  ergodic, almost mixing and mixing when the dominant eigenvalues of
  $f^*$ on $H^{p,p}(X,\R)$ are equal to $d_p$.}.
This measure is the intersection
of a Green current $T^+$ with $(1,1)$-currents with H{\"o}lder continuous potentials.
The main result in \cite{DinhNguyenSibony} implies that the measure is moderate.
Here is a criterion for the non-vanishing
of $T^+\wedge T^-$, see also \cite{Guedj}. 

\begin{proposition} \label{prop_green_current_wedge}
If $f$ is as above, then the following properties are
equivalent
\begin{enumerate}
\item There is a Green $(p,p)$-current $T^+$ of $f$ and a
      Green $(k-p,k-p)$-current $T^-$ of $f^{-1}$ such that
      $T^+\wedge T^-$ is a positive non-zero measure.
\item The spectral radius of $f^*$ on $H^{p,p}(X,\R)$ is of multiplicity $1$.
\item The spectral radius of $f_*$ on $H^{k-p,k-p}(X,\R)$ is of multiplicity $1$.
\end{enumerate}
\end{proposition}
\proof
To say that the multiplicity of the spectral radius is 1 means that
the Jordan blocks associated to the eigenvalues of maximal modulus are
reduced to these eigenvalues.
As it is showed in Section \ref{section_linear}, if we consider a basis of
$H^{p,p}(X,\R)$ and $H^{k-p,k-p}(X,\R)$ which are dual with respect to
the cup-product $\smile$, then $f^*$ acting on $H^{p,p}(X,\R)$ and $f_*$ acting on
$H^{k-p,k-p}(X,\R)$ are given by the same matrix. Therefore, 
properties 2 and 3 are equivalent.  

Assume that properties 2 and 3 hold. Then, by Proposition \ref{prop_cesaro_cv}, the currents
$$S_N^+:=\widehat L_N^+(\omega^p)\quad \mbox{where}\quad L_N^+:={1\over N}\sum_{n=1}^N d_p^{-n} (f^n)^*,$$
converge SP-uniformly to a positive closed $(p,p)$-current
$T^+$. By Lemma \ref{lemma_mass_current}, the cohomology class of
$T^+$ is non-zero. Therefore, $T^+$ is a Green $(p,p)$-current.
In the same way, we prove that 
$$S_N^-:=\widehat L_N^-(\omega^{k-p})\quad \mbox{where} \quad \widehat
L_N^-:={1\over N}\sum_{n=1}^N d_p^{-n} (f^n)_*,$$
converge SP-uniformly to a Green $(k-p,k-p)$-current $T^-$ of
$f^{-1}$. 

Since the convergences are SP-uniform, $S_N^+\wedge S_N^-$
converge to $T^+\wedge T^-$. Observe that $S_N^+\wedge S_N^-$
is a smooth positive measure.
In order to prove property
1, we only have to check that the mass of $S_N^+\wedge S_N^-$ does
not tend to 0. We have
\begin{eqnarray*}
\|S_N^+\wedge S_N^-\| & = & {1\over N^2}\sum_{1\leq n,l\leq N}
d_p^{-n-l}\int_X (f^n)^*(\omega^p)\wedge (f^l)_*(\omega^{k-p}) \\
& = & {1\over N^2}\sum_{1\leq n,l\leq N} d_p^{-n-l}\int_X (f^{n+l})^*(\omega^p)\wedge \omega^{k-p} \\
& = &  {1\over N^2}\sum_{1\leq n,l\leq N}
d_p^{-n-l}\|(f^{n+l})^*(\omega^p)\|.
\end{eqnarray*}
By Lemma \ref{lemma_mass_current}, the last quantity is bounded from below by a positive
constant independent of $N$. This implies that  the mass of $S_N^+\wedge S_N^-$ does
not tend to 0. 

Now assume property 1 and let $m$ denote the multiplicity of the
spectral radius of $f^*$ on $H^{p,p}(X,\R)$. By Proposition
\ref{prop_green_inverse},
there are smooth closed
$(p,p)$-form $S^+$ and $(k-p,k-p)$-form $S^-$, not necessarily
positive, and an increasing sequence $(n_i)$  such that
$$T^+=\lim_{i\rightarrow\infty} n_i^{1-m}d_p^{-n_i}
(f^{n_i})^*(S^+)\quad \mbox{and}\quad T^-=\lim_{i\rightarrow\infty} n_i^{1-m}d_p^{-n_i}
(f^{n_i})_*(S^-).$$
Moreover, the convergences are SP-uniform. We have
\begin{eqnarray*}
\|T^+\wedge T^-\| & = &  \lim_{i\rightarrow\infty}\int_X  n_i^{1-m}d_p^{-n_i}
(f^{n_i})^*(S^+)\wedge n_i^{1-m}d_p^{-n_i}
(f^{n_i})_*(S^-) \\
& = & \lim_{i\rightarrow\infty}\int_X  n_i^{2-2m}d_p^{-2n_i}
(f^{2n_i})^*(S^+)\wedge S^-\\
& = & \lim_{i\rightarrow\infty}\Big({2\over n_i}\Big)^{m-1}\int_X  (2n_i)^{1-m}d_p^{-2n_i}
(f^{2n_i})^*(S^+)\wedge S^-.
\end{eqnarray*}
The last integral can be computed cohomologically. By Lemma \ref{lemma_mass_current}, it
converges to a constant when $i$ tends to infinity. Therefore, if
$\|T^+\wedge T^-\|$ is strictly positive,
$(2/n_i)^{m-1}$ does not converge to 0. It follows that $m=1$. 
\endproof

When $p=1$, using Hodge-Riemann theorem one prove easily that the
properties in Proposition \ref{prop_green_current_wedge} are always satisfied.
We don't know if this is the case in general. For this question, the
reader will find some useful results
and techniques developed in \cite{Gromov1, DinhNguyen, KeumOguisoZhang}. 
Here is the main result of this section. The property that $\mu$ is of
maximal entropy was obtained in collaboration with de Th\'elin.

\begin{theorem} \label{th_measure}
Let $f$ be a holomorphic automorphism on a compact K{\"a}hler manifold
$(X,\omega)$ of dimension $k$ and $d_s$ the dynamical degrees
of $f$. Assume that there is a dynamical degree $d_p$ strictly
larger than the other ones and that $f$
satisfies the properties in Proposition
\ref{prop_green_current_wedge}. Then $f$ admits an invariant 
probability measure $\mu$ with H{\"o}lder continuous super-potentials.
The measure $\mu$ is ergodic, hyperbolic and of maximal entropy. 
If the
dominant eigenvalues of $f^*$ on $H^{p,p}(X,\C)$ are equal to $d_p$
then $\mu$ is mixing.
\end{theorem}

We will see that the last assertion holds under a weaker
hypothesis: {\it $d_p$ is the unique dominant eigenvalue
which is a root of a real number}. This condition is always satisfied
for some iterates of $f$. 
The measure $\mu$ is called {\it equilibrium measure} of $f$. By
Propositions \ref{prop_holder_moderate} and
\ref{prop_holder_hausdorff}, $\mu$ 
is moderate and has positive Hausdorff dimension when $X$ is
projective. Since $\mu$ is hyperbolic, a theorem by Katok
\cite[p.694]{KatokHasselblatt} 
says that any point in
the support of $\mu$ can be approximated by saddle periodic
points. Therefore, the saddle periodic points are Zariski dense
in $X$ since moderate measures have no mass on proper analytic subsets
of $X$. 

Recall that a positive invariant measure $\mu$ is {\it
  mixing} if for any test functions $\phi$, $\psi$ (smooth, continuous,
bounded or in $L^2(\mu)$) we have
$$\langle \mu, (\phi\circ f^n)\psi\rangle \rightarrow 
\|\mu\|^{-1}\langle\mu,\phi\rangle\langle\mu,\psi\rangle.$$ 
The invariance of $\mu$ implies that $\langle \mu, (\phi\circ f^n)\psi\rangle =\langle \mu, \phi(\psi\circ
f^{-n})\rangle$. So, $\mu$ is mixing for $f$ if and only if it is mixing for $f^{-1}$.
Mixing is equivalent to the the property that $(\phi\circ f^n)\mu$
converge to a constant times $\mu$. Indeed, since $\mu$ is invariant,
we have $\langle
(\phi\circ f^n)\mu,1\rangle =\langle \mu,\phi\rangle$; hence, the above
constant should be $\|\mu\|^{-1}\langle \mu,\phi\rangle$. 
In fact, mixing is also equivalent to the property that any limit value
of  $(\phi\circ f^n)\mu$ is proportional to $\mu$.

If $\mu$ is mixing then it is {\it ergodic}, that is,
$\mu$ is extremal in the cone of positive invariant measures. This
property is equivalent to the convergence
$${1\over N}\sum_{n=1}^N\langle \mu, (\phi\circ f^n)\psi\rangle \rightarrow 
\|\mu\|^{-1}\langle\mu,\phi\rangle\langle\mu,\psi\rangle$$
or to the property that any limit value of 
$$\mu_N:={1\over N}\sum_{n=1}^N (\phi\circ f^n) \mu$$
is proportional to $\mu$. Note that $\mu$ is ergodic for $f$ if and
only if it is
ergodic for $f^{-1}$.
We refer to \cite{KatokHasselblatt, Walters} for the notions of entropy and of Lyapounov
exponent. An invariant positive measure is {\it hyperbolic} if its
Lyapounov exponents are non-zero.  

We recall and introduce some notations that we will use. 
Let $F$ and $H$ be the dominant and strictly dominant subspaces of
$H^{p,p}(X,\R)$ for the action of $f^*$. Let $F^\vee$ and $H^\vee$ 
denote the dominant and strictly dominant subspaces of
$H^{k-p,k-p}(X,\R)$ for the action of $f_*$. By Theorem
\ref{th_green_construction}, we can 
associate to each class $c$ in $F$ or in $H$ a current in
$\Dc_p$ with H{\"o}lder continuous super-potentials that we denote by
$T_c^+$. We can also apply this result to $f^{-1}$ and associate to
each class $c^\vee$ in $F^\vee$ or in $H^\vee$ a current $T^-_{c^\vee}$
in $\Dc_{k-p}$ with H{\"o}lder continuous super-potentials. By
Proposition \ref{prop_holder_wedge}, the measure $T_c^+\wedge T_{c^\vee}^-$ has
H{\"o}lder continuous super-potentials. 

Denote by $\Mc$ the real space generated by $T_c^+\wedge
T_{c^\vee}^-$ with $c\in F$ and $c^\vee\in F^\vee$ and $\Nc$ the real space generated by $T_c^+\wedge
T_{c^\vee}^-$ with $c\in H$ and $c^\vee\in H^\vee$. We have
$$\dim \Mc\leq (\dim F)(\dim F^\vee)=(\dim F)^2$$
and
$$\dim \Nc\leq (\dim H)(\dim H^\vee)=(\dim H)^2.$$
Let $\Mc^+$, $\Nc^+$ be the closed convex cones of positive measures
in $\Mc$ and in $\Nc$. 
The measure $\mu$ that we will
construct is an extremal element of $\Nc^+$. 
We first prove the following lemmas.

\begin{lemma} \label{lemma_wedge_green_inv}
For all $c\in F$ and $c^\vee\in F^\vee$, we have 
$$f^*(T^+_c\wedge T^-_{c^\vee})=T^+_{f^*c}\wedge T^-_{f^*c^\vee}.$$
If $c$ is in $H$ and $c^\vee$ is in $H^\vee$, then $T_c^+\wedge
T_{c^\vee}^-$ is an invariant measure.
\end{lemma}
\proof
Write as in Proposition
\ref{prop_green_construction}
$$T_c^+=\lim_{i\rightarrow\infty} d_p^{-n_i}
(f^{n_i})^*(S^+)\quad \mbox{and}\quad T_{c^\vee}^-=\lim_{i\rightarrow\infty} d_p^{-n_i}
(f^{n_i})_*(S^-),$$
with $S^+$, $S^-$ smooth.
Observe that $d_p^{-n_i}
(f^{n_i})^*[S^+]$ converge to the class $c$ and  $d_p^{-n_i}
(f^{n_i})_*[S^-]$ converge to $c^\vee$. Hence, $d_p^{-n_i}
(f^{n_i})^*[f^*(S^+)]$ converge to $f^*c$ and
$d_p^{-n_i} (f^{n_i})_*[f^*(S^-)]$ converge to
$f^*c^\vee$. 
Applying Proposition \ref{prop_green_construction} to $f^*(S^+)$ and to $f^*(S^-)$ yields
\begin{eqnarray*}
f^*(T^+_c\wedge T^-_{c^\vee}) & = & f^*\big(\lim_{i\rightarrow\infty} d_p^{-n_i}
(f^{n_i})^*(S^+)\wedge d_p^{-n_i} (f^{n_i})_*(S^-)\big)\\
& = & \lim_{i\rightarrow\infty} d_p^{-n_i}
(f^{n_i})^*\big(f^*(S^+)\big)\wedge d_p^{-n_i}
(f^{n_i})_*\big(f^*(S^-)\big)\\
& = & T^+_{f^*c}\wedge T^-_{f^*c^\vee}.
\end{eqnarray*}
When $c$ is in $H$ and $c^\vee$ is in $H^\vee$, we have $f^*c=d_pc$
and $f^*c^\vee= d_p^{-1} c^\vee$. Therefore, $T^+_{f^*c}=d_pT_c^+$ and
$T^-_{f^*c^\vee}=d_p^{-1} T^-_{c^\vee}$. We deduce that
$f^*(T^+_c\wedge T^-_{c^\vee})=T^+_c\wedge T^-_{c^\vee}$. 
Since $f$
is an automorphism, this also implies that $f_*(T^+_c\wedge T^-_{c^\vee})=T^+_c\wedge T^-_{c^\vee}$.
Hence, $T^+_c\wedge T^-_{c^\vee}$ is invariant.
\endproof

\begin{lemma}
The cones $\Mc^+$ and $\Nc^+$ have non-empty interior in $\Mc$ and $\Nc$
respectively.
\end{lemma}
\proof
Consider $c\in F$ and $c^\vee\in F^\vee$. 
Observe that if $T_c^+$ or $T_{c^\vee}^-$ is approximable by smooth positive closed
currents then $T_c^+\wedge T_{c^\vee}^-$ is positive and belongs to
$\Mc^+$. 
We have seen that the sets of such classes have non-empty
interiors in $F$ and $F^\vee$.  Hence, $\Mc$ is generated by such positive measures $T_c^+\wedge T_{c^\vee}^-$.
It follows that $\Mc^+$ has
non-empty interior. The case of $\Nc^+$ is treated in the same way.
\endproof

\begin{definition} \rm Let $\mu$ be an invariant positive measure of
  $f$. We say that $\mu$ is {\it almost mixing} if
there is a finite dimensional space $V$ of measures such that for any
  function $\phi$ in $L^2(\mu)$ the limit values of $(\phi\circ f^n)\mu$, when $n\rightarrow\infty$,
  belong to $V$.
\end{definition}
The above notion does not change if we use the
  continuous functions or the space $L^1(\mu)$ instead of $L^2(\mu)$
  since continuous and $L^2$ functions are dense in $L^1(\mu)$.
Note also that mixing corresponds to the case where $V$ is of
dimension 1. We will see in the following lemma that $\mu$ is almost mixing if and
only if $L^2(\mu)$ can be decomposed into an invariant orthogonal sum
$W\oplus W^\perp$ with $W^\perp$ of finite dimension such that
$(\phi\circ f^n)\mu\rightarrow 0$ for $\phi\in W$. We can deduce that
$(\psi\circ f^{-n})\mu\rightarrow 0$ for $\psi\in W$ and that
$\mu$ is also almost mixing for $f^{-1}$.
The following lemma is valid for a general dynamical system.

\begin{lemma} \label{lemma_ergodic_mixing}
Let $\mu$ be a positive measure invariant by
  $f$. Assume that $\mu$ is almost mixing and that 
 $\mu$ is ergodic for every $f^n$ with $n\geq 1$. Then $\mu$ is mixing.  
\end{lemma}
\proof
Let $V$ be the smallest space 
of measures such that for any
  real-valued continuous
  function $\phi$ the limit values of $(\phi\circ f^n)\mu$, when $n\rightarrow\infty$,
  belong to $V$. This space is invariant by $f^*$ and $f_*$.
We have to prove that $\dim V=1$.
Let $W$ denote the space of functions $\psi\in L^2(\mu)$ with complex
values such that 
$\langle \mu',\psi\rangle=0$ for every $\mu'\in V$. Let $W^\perp$
denote the orthogonal of $W$.
The spaces $W$, $W^\perp$ are invariant under $f^*$, $f_*$ and
we have $\dim_\C W^\perp =\dim V$. 
Moreover, continuous functions are dense in $W$.
We show that $\dim_\C W^\perp=1$.

Since $f^*$ and $f_*$ preserve the scalar product in $L^2(\mu)$, all the eigenvalues 
of $f^*$ and of $f_*$ have modulus equal to 1. So, if $\psi$ is an eigenvector
of $f^*$ associated to an eigenvalue $\lambda$, we have
$|\psi|\circ f=|\psi|$ and then $|\psi|$ is constant since $\mu$ is ergodic.
Therefore, $\psi^n\in L^2(\mu)$ and 
$\psi^n\circ f=\lambda^n\psi^n$ for every $n\in\Z$. 
We claim that $W$ does not contain any eigenvector. Otherwise, there is 
a function $\psi\in W\setminus\{0\}$ such that $\psi\circ
f=\lambda\psi$ with 
$|\lambda|=1$. 
Since $\psi$ can be approximated by continuous functions in $W$, 
we deduce from the definition of $W$ that for every $\phi\in L^2(\mu)$:
$$|\langle \psi\mu,\phi\rangle|=|\langle (\psi\circ f^{-n})\mu,\phi\rangle|
=|\langle (\phi\circ f^n)\mu,\psi\rangle|\rightarrow 0.$$
We get that $\psi\mu=0$, hence $\psi=0$. This is a contradiction.

Consider now an eigenvector $\psi$ of $f^*$ in $W^\perp$ associated to an eigenvalue
$\lambda$. Then, $\psi^n$ is an eigenvector associated 
to $\lambda^n$ for every $n\in \Z$. We deduce that $\psi^n$ is a
function in $W^\perp$.  Since $W^\perp$ is finite dimensional, $\lambda$ 
is a root of unity. We have $\psi\circ f^n=\psi$ for some $n\geq 1$.
Since $\mu$ is ergodic for $f^n$, $\psi$ is constant. Hence, $\lambda=1$ 
and it follows that $\dim_\C W^\perp=1$ because $f^*$ is an isometry of $W^\perp$. 
\endproof

Consider the automorphism $\widetilde f$ of $X\times X$ given by
$\widetilde f(x,y):=(f(x),f^{-1}(y))$. By K{\"u}nneth formula, we have
$$H^{l,l}(X\times X,\C)\simeq \sum_{r+s=l} H^{r,s}(X,\C)\otimes
H^{s,r}(X,\C).$$
Moreover,  $\widetilde f^*\simeq (f^*,f_*)$ preserves this
decomposition. It is shown in \cite{Dinh1} that the spectral radius of
$f^*$ on $H^{r,s}(X,\C)$ is bounded by $\sqrt{d_rd_s}$. We deduce that
$d_k(\widetilde f)$ is the maximal dynamical degree of $\widetilde
f$. It is equal to $d_p^2$, with multiplicity 1 and is strictly larger
than the other ones. So, the results obtained for $f$ can be applied
to $\widetilde f$. We will deduce several properties for $f$. 

We 
use analogous notations $\widetilde \Nc$,
$\widetilde \Nc^+$... for $\widetilde f$ instead of the notations $\Nc$,
$\Nc^+$... for $f$. By Theorem \ref{th_green_construction} and
Corollary \ref{cor_unique_green} applied to $\widetilde f$, together with the
K{\"u}nneth formula, the family of the Green $(k,k)$-currents of
$\widetilde f$ is a convex cone with non-empty interior in the real
space generated by the currents $T_c^+\otimes T_{c^\vee}^-$. The Green $(k,k)$-currents of
$\widetilde f^{-1}$ is a convex cone with non-empty interior in the real
space generated by the currents $T_{c^\vee}^-\otimes
T_c^+$. Therefore, $\widetilde \Nc$ is generated by $\mu\otimes\mu'$
with $\mu$, $\mu'$ in $\Nc$ and $\widetilde \Mc$ is generated by
$\mu\otimes \mu'$ with $\mu$, $\mu'$ in $\Mc$.

Let $S^+$
be a smooth current in $\Dc_p$ such that $d_p^{-n_i}(f^{n_i})^*S^+$
converge SP-uniformly to $T^+_c$ for some increasing sequence
$(n_i)$. Let $S^-$ be a smooth current in $\Dc_{k-p}$ such that $d_p^{-n_i}(f^{n_i})_*S^-$
converge SP-uniformly to $T^-_{c^\vee}$. Then, we deduce from Proposition
\ref{prop_green_construction} applied to $\widetilde f$ 
that $d_p^{-2n_i}(\widetilde f^{n_i})^*(S^+\otimes S^-)$
converge SP-uniformly to $T_c^+\otimes T_{c^\vee}^-$. We will use this
property in the computations involving  $T_c^+\otimes
T_{c^\vee}^-$.

\begin{lemma} \label{lemma_stable_N}
Let $\phi$ be a continuous real-valued
function on $X$. 
If $\mu$ is a measure in $\Mc$, then any limit value of $(\phi\circ
f^n)\mu$ is a measure in $\Mc$. In particular, the measures in $\Nc^+$
are almost mixing. If $\mu$ is in $\Nc$, then
any limit value of 
$$\mu_N:={1\over N}\sum_{n=1}^N (\phi\circ f^n)\mu$$
is a measure in $\Nc$.
\end{lemma}
\proof
Since continuous functions are uniformly approximable by smooth
functions, we can assume that $\phi$ is smooth.
We prove the first assertion. By definition of $\Mc$, we can assume
that $\mu=T^+\wedge T^-$ where $T^+$ is a $(p,p)$-current associated
to a class $c$ in $F$
and $T^-$ is a $(k-p,k-p)$-current associated to a class $c^\vee$ in
$F^\vee$ as above.
It is enough to show that if a subsequence 
$(\phi\circ f^{2n_i})\mu$ converge, then the limit is a
measure in $\Mc$. Indeed, we obtain the case with odd powers by
replacing $\phi$ by
$\phi\circ f$.

Let $\psi$ be another test smooth function on $X$. Define
$\Phi(x,y):=\phi(x)\psi(y)$. Since $\mu$ is
invariant, lifting the integrals on $X$ to $\Delta$ we get
\begin{eqnarray*}
\langle (\phi\circ f^{2n})\mu,\psi\rangle & = & \langle \mu,
(\phi\circ f^n)(\psi\circ f^{-n}) \rangle \\
& = & \langle T^+\wedge T^-,(\phi\circ f^n)(\psi\circ f^{-n})
\rangle\\
& = & \langle (T^+\otimes T^-)\wedge [\Delta],\Phi\circ \widetilde f^n\rangle,
\end{eqnarray*}
where in order to obtain the last line we use a SP-uniform approximation of $T^+\otimes
T^-$ by smooth currents as above. We have
\begin{eqnarray*}
\langle (\phi\circ f^{2n})\mu,\psi\rangle & = &  \langle (\widetilde
f^n)_*(T^+\otimes T^-)\wedge (\widetilde f^n)_*[\Delta],\Phi\rangle \\
& = &  \big\langle (d_p^n(f^n)_*T^+\otimes d_p^n(f^n)^*T^-)\wedge
d_p^{-2n} (\widetilde f^n)_*[\Delta],\Phi\big \rangle.
\end{eqnarray*}
Observe that $d_p^n(f^n)_*T^+$ belongs to a bounded family of
currents constructed in Theorem \ref{th_green_construction}. An analogous property holds
for $d_p^n(f^n)^*T^-$. Therefore, 
the limit values of
$$(d_p^n(f^n)_*T^+\otimes d_p^n(f^n)^*T^-)\wedge d_p^{-2n} (\widetilde f^n)_*[\Delta]$$
are measures in
$\widetilde\Mc$. It follows that $\langle (\phi\circ
f^{2n_i})\mu,\psi\rangle$ converge to a finite combination of 
$$\langle \mu^+\otimes\mu^-,\Phi\rangle = \const
\langle\mu^-,\psi\rangle$$
with $\mu^+$, $\mu^-$ in $\Mc$. We deduce that $(\phi\circ f^{2n_i})\mu$
converge to a combination of $\mu^-$. So, the limit values of
$(\phi\circ f^{2n})\mu$ are in $\Mc$. This
completes the proof of the first assertion.

For the last assertion, we follow the same approach with $T^+$
associated to a class in $H$ and $T^-$ associated to a class in
$H^\vee$. In this case, $T^+$, $T^-$ are invariant and any limit value of 
$$(T^+\otimes T^-)\wedge {1\over N} \sum_{n=1}^Nd_p^{-2n}
(\widetilde f^n)_*[\Delta]$$
is a measure in $\widetilde \Nc$. We deduce as above that the
limit values of $\mu_N$ are in $\Nc$.
\endproof

\begin{proposition} \label{prop_mixing}
Let $\mu$ be a probability measure in $\Nc^+$. Then $\mu$ is ergodic if
and only if it is an extremal element of $\Nc^+$. 
Moreover, the number of extremal probability measures in $\Nc^+$ is equal to $\dim
\Nc$ and the convex cone $\Nc^+$ is generated by these measures.
When $d_p$ is the only dominant
eigenvalue of $f^*$ on $H^{p,p}(X,\C)$ which is a root of a real
number, then $\mu$ is mixing if
and only if it is an extremal element of $\Nc^+$.
\end{proposition}
\proof
If $\mu$  ergodic, $\mu$ is extremal in the cone of invariant positive
measures. Therefore, $\mu$ is extremal in
$\Nc^+$. Assume now that $\mu$ is extremal in $\Nc^+$. We show that it
is ergodic.
Let $\phi$ be a positive continuous function. 
The measures $\mu_N$, defined as above, are positive and bounded by $\|\phi\|_\infty\mu$.
By Lemma \ref{lemma_stable_N}, any
limit value of $\mu_N$
is a measure in $\Nc^+$ and it is bounded by  $\|\phi\|_\infty\mu$.
Since $\mu$ is extremal in $\Nc^+$, these limit
values are proportional to $\mu$. 
Therefore, $\mu$ is ergodic. 

Recall that $\Nc^+$ is a salient convex closed cone in
$\Nc$ with non-empty interior. Moreover, 
any element $\nu$ of $\Nc^+$ is an integral over extremal elements of
mass 1. So, we get a decomposition of $\nu$ into ergodic probability
measures. Since this decomposition is unique \cite{Walters} and since $\Nc$ is
generated by $\dim\Nc$ elements,
we deduce that the
number of extremal probability measures in $\Nc^+$ is equal to $\dim
\Nc$ and the convex cone $\Nc^+$ is generated by these measures. So,
$\Nc^+$ is a cone with simplicial basis.  

Assume that $d_p$ is the only dominant
eigenvalue of $f^*$ on $H^{p,p}(X,\C)$ which is a root of a real
number. Then the spaces $H$, $H^\vee$ do not change if we replace
$f$ by $f^n$. Therefore, $\Nc$ do not change if we replace
$f$ by $f^n$.
We deduce that $\mu$ is ergodic for $f^n$.
Lemmas \ref{lemma_stable_N} and
\ref{lemma_ergodic_mixing} imply that $\mu$ is mixing. This completes
the proof of the proposition.
\endproof

\noindent
{\bf End of the proof of Theorem \ref{th_measure}.} Let $\mu$ be a
probability measure which is an
extremal element of $\Nc^+$. By Proposition \ref{prop_mixing}, $\mu$ is ergodic and
is mixing if $d_p$ is the only dominant eigenvalue of $f^*$ on
$H^{p,p}(X,\C)$ which is a root of a real number. By definition of
$\Nc$, this measure has H\"older continuous super-potentials. It
remains to prove that $\mu$ is of maximal entropy. Indeed, by a recent
result of de Th\'elin \cite{deThelin}, the property that $\mu$ is of
entropy $\log d_p$  together with the fact that $d_p$ is
strictly larger than the other dynamical degrees implies that $\mu$ is
hyperbolic, see also \cite{DinhNguyenSibony1}. More precisely, $\mu$ admits $p$ positive Lyapounov exponents
larger than or equal to ${1\over 2}\log (d_p/d_{p-1})$ and $k-p$ negative
  exponents at most equal to $-{1\over 2}\log(d_p/d_{p+1})$. 

The variational principle
\cite{Walters} implies that the entropy of an invariant measure is
bounded from above by the topological entropy of $f$. 
By Gromov and Yomdin results \cite{Gromov2,Yomdin}, the topological
entropy of $f$ is equal to $\log d_p$. Therefore, if $\nu$ is a probability
measure in $\Nc^+$ then the entropy $h(\nu)$ of $\nu$ is at most equal
to $\log d_p$. We will prove that $h(\nu)=\log d_p$ for every
probability measure $\nu$ in
$\Nc^+$.

Let $S^+$ be a smooth form in $\Dc_p$ and $S^-$
a smooth form in $\Dc_{k-p}$. If $S^+$ and $S^-$ are strictly positive, 
by Proposition \ref{prop_entropy_dt} in the appendix below, any limit value $\nu$ of 
$$\nu_n:={1\over n}\sum_{l=1}^n d_p^{-n} (f^l)^*(S^+)\wedge
(f^{n-l})_*(S^-)$$
is proportional to an invariant probability measure of maximal
entropy $\log d_p$. By Proposition \ref{prop_green_construction}, 
the space $M$ generated by these measures
$\nu$ is of finite dimension. Let $M_P$ denote the convex of probability
measures in $M$.
Since the entropy $h(\nu)$ is an affine  function on $\nu$
\cite[p.183]{Walters}, all the measures in $M_P$ are of entropy $\log d_p$. It
suffices to show that $M$ contains $\Nc$. Observe that since $\nu_n$
depends linearly on $S^+$, $S^-$, the space $M$ contains
also the limit values of $\nu_n$ when $S^+$, $S^-$ are not necessarily
positive. When $S^+$ is in a class $c\in H$ and $S^-$ is in a class
$c^\vee\in H^\vee$, by Proposition \ref{prop_green_construction},
$\nu_n$ converge to $T_c^+\wedge T_{c^\vee}^-$. We deduce that $M$
contains $\Nc$ and this implies the result.
\hfill $\square$

\begin{remark} \rm
The property that the equilibrium measures are of maximal entropy can
be proved using $\widetilde f$. More precisely, using
Proposition \ref{prop_entropy_bs} below for $Y:=\Delta$, we can
construct equilibrium measures of $\widetilde f$ with maximal entropy
$2\log d_p$. This together with the Brin-Katok formula applied to $f$,
$f^{-1}$ and $\widetilde f$, see the appendix below and \cite[p.99]{Walters}, implies that the
equilibrium measures of $f$ are of entropy $\log d_p$.  The use
of $\widetilde f$ may be a good method in order to study the
distribution of periodic points of $f$ by considering the
intersection $(\widetilde f^n)^* [\Delta]\wedge [\Delta]$.
\end{remark}

\begin{remark}\rm
The Green currents and the equilibrium measures have been constructed and
studied by the authors in \cite{DinhSibony4}, for $f$ with a
dynamical degree $d_p$ strictly larger than the other ones. Guedj
considered in \cite{Guedj} the situation with the aditional hypothesis that 
$d_p$ is the unique dominant eigenvalue of $f^*$ on
$H^{p,p}(X,\C)$, i.e. $\dim F=\dim H=1$. He claims that when $X$ is
projective, the equilibrium measure is of maximal entropy but he didn't give
the proof. In this situation, we can find a subvariety $Y$ of
dimension $p$ in $X$ such that $d_p^{-n}(f^n)_*[Y]$ converge to a
Green current $T^-$, see also \cite{DinhSibony6}. Then, using the
SP-uniform convergence $d_p^{-n}(f^n)^*(\omega^p)\rightarrow
  T^+$ or properties proved in \cite{DinhSibony4}, we deduce that 
$d_p^{-n-l}(f^n)^*(\omega^p)\wedge
(f^l)_*[Y]$ converge to a constant times the equilibrium measure which, by
Proposition \ref{prop_entropy_bs} below, is of maximal entropy. 
\end{remark}


\subsection{Appendix: measures of maximal entropy} \label{section_appendix}

This section contains an abstract construction of measures of maximal
entropy. Most of the arguments given here are well-known, see
Bedford-Smillie \cite{BedfordSmillie} and de Th\'elin
\cite{deThelin2}. For simplicity, assume that $f:X\rightarrow X$ 
is an automorphism as above which satisfies the properties in
Proposition \ref{prop_green_current_wedge}. The last hypothesis garantees
that the construction gives non-zero measures. The method is
still valid in a much more general setting, in particular, when $f$
is a non-invertible finite map.  

Given $\epsilon>0$ and $n\in\N$, define {\it the Bowen ball}
$B_n(a,\epsilon)$ by 
$$B_n(a,\epsilon):=\big\{x\in X,\quad \dist(f^i(x),f^i(a))\leq
\epsilon \quad \mbox{for}\quad 0\leq i\leq n\big\}.$$
Let $\nu$ be a probability measure invariant by $f$. 
By Brin-Katok \cite{BrinKatok}, the function
$$h(\nu,a):=\sup_{\epsilon>0}\lim_{n\rightarrow\infty}-{1\over
  n}\log\nu(B_n(a,\epsilon))$$
is well-defined $\nu$-almost everywhere and the entropy of $\nu$ is equal to
$$h(\nu)=\int h(\nu,a) d\nu(a).$$
We have the following Misiurewicz's lemma which is valid for continuous
maps on compact metric spaces, see \cite{BedfordSmillie,deThelin2,Walters}.

\begin{lemma} \label{lemma_misiu}
Let $(n_i)$ be an increasing sequence of integers and $\nu_{n_i}$
probability measures such that the sequence
$${1\over n_i}\sum_{l=0}^{n_i-1} (f^l)_*(\nu_{n_i})$$
converges to a measure $\nu$. Assume there are constants $\epsilon>0$
and $c_{n_i}>0$ such that 
$\nu_{n_i}(B_{n_i}(a,\epsilon))\leq c_{n_i}$ for all $i$ and all Bowen
ball $B_{n_i}(a,\epsilon)$.
Then $\nu$ is an invariant probability
measure and its entropy $h(\nu)$ satisfies the inequality
$$h(\nu)\geq \limsup_{i\rightarrow\infty}-{1\over n_i}\log c_{n_i}.$$
\end{lemma}

We deduce from this lemma and an estimate due to Yomdin \cite{Yomdin}
the following proposition which was obtained in collaboration with de
Th\'elin.

\begin{proposition} \label{prop_entropy_dt}
Let $S^+$ be a bounded positive $(p,p)$-form and $S^-$ a bounded
positive $(k-p,k-p)$-form on $X$, not necessarily closed. Assume
there is an increasing sequence $(n_i)$ of integers such that
$${1\over n_i}\sum_{l=1}^{n_i} d_p^{-n_i} (f^l)^*(S^+)\wedge
(f^{n_i-l})_*(S^-)$$
converge to a probability measure $\nu$. Then $\nu$ is an invariant
measure of maximal entropy $\log d_p$.
\end{proposition}
\proof
Denote by $\nu_{n_i}'$ the positive measure $d_p^{-n_i}
(f^{n_i})^*(S^+)\wedge S^-$. Define 
$\nu_{n_i}:=\lambda_{n_i}^{-1}\nu_{n_i}'$ where
$\lambda_{n_i}$ is the mass of $\nu_{n_i}'$. Then $\nu_{n_i}$ are
probability measures and we have 
$$\lambda_{n_i}{1\over n_i}\sum_{l=0}^{n_i-1} (f^l)_*(\nu_{n_i})=
{1\over n_i}\sum_{l=1}^{n_i} d_p^{-n_i} (f^l)^*(S^+)\wedge
(f^{n_i-l})_*(S^-)$$
which converge to the probability measure $\nu$. We deduce that
$\lambda_{n_i}$ converge to 1.
Therefore, by Lemma \ref{lemma_misiu}, it is enough to prove for
any $0<\delta<1$
the existence of positive constants $\epsilon$, $A$ such that
$\nu_{n_i}'(B_{n_i}(a,\epsilon))\leq A d_p^{-n_i}e^{\delta n_i}$ for every
$a\in X$. 
For this purpose, we can assume for simplicity that $S^+=\omega^p$ and
$S^-=\omega^{k-p}$.  We have to show that $\nu_n''(B_n(a,\epsilon))\leq A
e^{n\delta}$ where $\nu_n'':=(f^n)^*(\omega^p)\wedge
\omega^{k-p}$. This inequality will be obtained by taking an average
on an estimate due to
Yomdin.

Let $Y\subset X$ be a complex manifold of dimension $p$ smooth up to
the boundary. If $\nu_n^Y:= (f^n)^*(\omega^p)\wedge [Y]$ then
 $\nu_n^Y(B_n(a,\epsilon))$ is equal to the volume of
 $f^n(Y\cap B_n(a,\epsilon))$ counted with multiplicity. Yomdin proved in \cite{Yomdin}
 that this volume is bounded by $Ae^{n\delta}$ when $\epsilon$ is
 small and $A$ is large enough. The estimate is uniform on $a$ and on
 $Y$. Now, consider a coordinate system $x=(x_1,\ldots,x_k)$ on a
 fixed chart of $X$ with $|x_i|<2$. In the unit polydisc $D$, up to a
 multiplicative constant, $\omega^{k-p}$ is bounded by
 $(\ddc\|x\|^2)^{k-p}$ which is equal to a combination of 
$$(idz_{i_1}\wedge d \overline z_{i_1})\wedge\ldots\wedge
(idz_{i_{k-p}}\wedge d \overline z_{i_{k-p}})\quad \mbox{with } 1\leq
i_1<\cdots<i_{k-p}\leq k.$$ 
The last form is equal to an average on the currents of
integration on the complex
submanifolds of $D$ which are given by
$$x_{i_1}=a_1,\quad\ldots, \quad x_{i_{k-p}}=a_{k-p} \quad \mbox{with
} a_i\in\C.$$
So, by Yomdin's inequality, $\nu_n''$ restricted to $D$ satisfies $\nu''_{n|D}(B_n(a,\epsilon))\leq
Ae^{n\delta}$ for some constants $\epsilon$, $A$. Since $X$
can be covered by a finite family of open sets $D$, we deduce that $\nu''_n(B_n(a,\epsilon))\leq
Ae^{n\delta}$ with $A>0$. This completes the proof. 
\endproof

One can prove in the same way the following proposition which is
essentially due to Bedford-Smillie \cite{BedfordSmillie}.

\begin{proposition} \label{prop_entropy_bs}
Let $S$ be a continuous positive $(p,p)$-form, $Y$ a complex
manifold of dimension $p$ in $X$ smooth up to the boundary and $\chi$
a bounded positive function on $Y$. Assume
there is an increasing sequence $(n_i)$ such that
$${1\over n_i}\sum_{l=1}^{n_i} d_p^{-n_i} (f^l)^*(S)\wedge
(f^{n_i-l})_*(\chi[Y])$$
converge to a probability measure $\nu$. Then $\nu$ is an invariant
measure of maximal entropy $\log d_p$.
\end{proposition}

More general situations will be considered by de Th\'elin and Vigny in
a forthcoming paper.


T.-C. Dinh, UPMC Univ Paris 06, UMR 7586, Institut de
Math{\'e}matiques de Jussieu, F-75005 Paris, France. {\tt
  dinh@math.jussieu.fr}, {\tt http://www.math.jussieu.fr/$\sim$dinh} 

\

\noindent
N. Sibony,
Universit{\'e} Paris-Sud, Math{\'e}matique - B{\^a}timent 425, 91405
Orsay, France. {\tt nessim.sibony@math.u-psud.fr}

\end{document}